\documentclass[10pt,aps,prb, preprint,reqno]{amsart}
\usepackage{amsmath}
\usepackage{amssymb}

\newtheorem{theorem}{Theorem}[section]
\newtheorem{remark}{Remark}[section]
\newtheorem{lemma}[theorem]{Lemma}
\newtheorem{proposition}[theorem]{Proposition}

\begin{document}
\title[Two dimensional MHD-alpha system]{A remark on the two-dimensional magneto-hydrodynamics-alpha system}
\author{Kazuo Yamazaki}  
\date{}
\thanks{The author expresses gratitude to Professor Jiahong Wu and Professor David Ullrich for their teaching, and the referee for helpful comments that helped the manuscript greatly.}
\maketitle

\begin{abstract}
We study the two-dimensional generalized magneto-hydrodynamics-$\alpha$ system with fractional Laplacians in the dissipative and diffusive terms. We show that the solution pair of velocity and magnetic fields preserves their initial regularity in all cases when the powers add up to one. This settles the global regularity issue in the general case which was remarked by the authors in [33] to be a problem. 

\vspace{5mm}

\textbf{Keywords: Magneto-hydrodynamics system, Navier-Stokes equations, global regularity, fractional Laplacians, $\alpha$-regularization}
\end{abstract}
\footnote{2010MSC : 35B65, 35Q35, 35Q86}
\footnote{Department of Mathematics and Statistics, Washington State University, Pullman, WA, 99164-3113, U.S.A.}

\section{Introduction}

We study the following two-dimensional generalized magneto-hydrodynamics (MHD)-$\alpha$ system: 
\begin{subequations}
\begin{align}
&\frac{\partial v}{\partial t} + (u\cdot\nabla) v + \sum_{k=1}^{2} v_{k}\nabla u_{k} + \nabla \left(\pi + \frac{1}{2} \lvert b\rvert^{2}\right) + \nu \Lambda^{2r_{1}}v = (b\cdot\nabla) b,\\
&\frac{\partial b}{\partial t} + (u\cdot\nabla) b - (b\cdot\nabla) u + \eta \Lambda^{2r_{2}} b = 0, \\
&v = (1-\alpha^{2}\Delta) u, \hspace{3mm} \nabla\cdot u = \nabla \cdot b = 0, \hspace{3mm} (v,b)(x,0) = (v_{0}, b_{0})(x),
\end{align}
\end{subequations}
where we denote $v = (v_{1}, v_{2}), u = (u_{1}, u_{2}), b = (b_{1}, b_{2})$ and $\pi$, the velocity, the filtered velocity, the magnetic and the pressure fields respectively. We also denote $\nu, \eta \geq 0$ the viscosity and diffusivity constants respectively and the fractional Laplacians defined through Fourier transform by 
\begin{equation*}
\widehat{\Lambda^{r_{i}}f}(\xi) = \lvert \xi\rvert^{r_{i}} \hat{f}(\xi), \hspace{3mm} i = 1, 2.
\end{equation*}
with their powers $r_{i} \geq 0$. Hereafter for simplicity let us write $\partial_{t} = \frac{\partial}{\partial t}, \partial_{i} = \frac{\partial}{\partial x_{i}}, i = 1, 2$. Finally, $\alpha > 0$ is the length-scale parameter representing the width of the filters. 

Let us briefly discuss the rich history concerning the MHD$-\alpha$ system. Firstly, when $b \equiv 0, r_{1} = 1$, the system (1a)-(1c) reduces to the viscous Camassa-Holm equations introduced in [5] which is well-known for its remarkable performance as a closure model of turbulence in infinite channels and pipes whose solutions give excellent agreement with empirical data for a wide range of large Reynolds numbers (cf. also [34]). It has much connections with the Navier-Stokes equations (NSE) of which its significance in engineering applications and mathematical difficulty in showing the global regularity result is well-known. 

Secondly, the system (1a)-(1c) is closely related to the generalized MHD system studied intensively since the work in [26]. The MHD system describes the motion of electrically conducting fluids and has broad applications in applied sciences such as astrophysics, geophysics and plasma physics (cf. [23]). The mathematical analysis of the MHD system has attracted much attention in particular regarding the global regularity issue in two dimension (cf. [1]-[2], [14]-[15], [24]- [25], [27]-[28], [30]-[32]). 

Finally, the MHD$-\alpha$ system (1a)-(1c) of our main concern was introduced in [20] with $r_{i} = 1, i = 1, 2$, purposely not filtering the magnetic field in contrast to the previous work. In a three-dimensional periodic domain, the authors in [20] obtained the global existence of the unique weak solution pair $(v, b)$ via Galerkin approximation and \textit{a priori } estimates and also obtained convergence as $\alpha \to 0$. In fact, it is shown in [36] that when $\nu, \eta > 0, r_{i} = 1 \hspace{1mm} \forall i$, the MHD$-\alpha$ system allows a unique weak solution to exist even in dimension four. We also refer to [6], [10], [13] and references therein for relevant work on related systems such as the Leray$-\alpha$ model, modified Leray-$\alpha$ subgrid scale model of turbulence and Bardina models. In particular that in [22] fractional Laplacians with various powers were taken into account in the study of such regularized models in a three-dimensional torus. Moreover, these regularized models have had much influence in the study of related equations, e.g. active scalars including the surface quasi-geostrophic equations (e.g. [7], [19], [29]). 

We now motivate our study. In [9] the authors showed that for the system (1a)-(1c), given $(v_{0}, b_{0}) \in H^{3}(\mathbb{R}^{2})$, the solution remains in such a space for all time in two cases: $\nu > 0, \eta = 0, r_{1} = 1  \hspace{3mm} \text{ or } \hspace{3mm} \nu = 0, \eta > 0, r_{2} = 1$. Similar results on a closely related system is also shown in [35]. In [33] the authors in particular observed that when $\nu, \eta > 0, r_{1} = r_{2} = \frac{1}{2}$, the solution pair to the system (1a)-(1c) remains smooth for all time. However, it was stated in [33] that the global regularity issue in case $r_{1} + r_{2} = 1$ in general is a  problem. We give an affirmative solution to this problem: 

\begin{theorem}
Suppose $\nu, \eta > 0, r_{1}, r_{2} \in (0, 1)$ so that $r_{1} + r_{2} = 1$. Then for any $(v_{0}, b_{0}) \in H^{3}(\mathbb{R}^{2})$, the solution pair $(v, b)$ to (1a)-(1c) remains smooth for all time: 
\begin{align*}
&v \in L^{\infty}([0, \infty); H^{3}(\mathbb{R}^{2})) \cap L^{2}([0, \infty); H^{3+r_{1}}(\mathbb{R}^{2})),\\
&b \in L^{\infty}([0,\infty ); H^{3}(\mathbb{R}^{2})) \cap L^{2}([0, \infty); H^{3+r_{2}}(\mathbb{R}^{2})).
\end{align*}
\end{theorem}

\begin{remark}
\begin{enumerate}
\item Theorem 1.1 solves the problem remarked in [33] Remark 1.1 and extends the result of [9] as well. Moreover, because the system (1a)-(1c) at $b \equiv 0, r_{1} = 1$ reduces to the viscous Camassa-Holm equation, our results also extend the study in that direction of research. 
\item As can be seen in the subsequent sections, the proof of Theorem 1.1 indeed required a series of delicate estimates. The case $r_{1} > \frac{1}{2}$ displayed difficulty due to the fact that because the velocity is filtered and dissipation strength is relatively stronger than that of diffusion, one is inclined to estimate the velocity equation first before the magnetic field equation. Indeed, one of the special feature of the system (1a)-(1c) is the simplicity of the vorticity formulation (see (21)) so that in its $L^{2}$-estimate, the only non-linear term that needs to be dealt with is $\nabla\times \left((b\cdot\nabla)b\right)$. The challenge is how to raise the regularity of $b$ to a level that we can handle this term, in particular when the diffusivity strength is relatively weak. Even if one compromises and estimates the velocity field before vorticity, the term $(b\cdot\nabla)b$ gives a problem. In the case $r_{2} > \frac{1}{2}$, it was also crucial to rely on the vorticity formulation (see Proposition 4.2). 
\item It seems to be an interesting problem whether the global regularity result can be extended below the threshold of $r_{1} + r_{2} = 1$. 
\item We would like to note the work on the two-dimensional Boussinesq system in [3], [11]-[12] by which our work was partially inspired. 
\end{enumerate}
\end{remark}

In the Preliminaries section, let us set up notations, state some key facts and useful lemmas. Subsequently, we prove \textit{a priori } estimates in the case $\frac{1}{2} < r_{1} < 1, 0 < r_{2} < \frac{1}{2}$ and then the case $0 < r_{1} < \frac{1}{2}, \frac{1}{2} < r_{2} < 1$. The case $r_{1} = r_{2} = \frac{1}{2}$ is done in [33]. Thereafter we conclude discussing the proof of Theorem 1.1. 

\section{Preliminaries}

Without loss of generality, let us assume $\nu = \eta = \alpha = 1$. Let us use the notation $A\lesssim_{a, b}B, A \approx_{a,b} B$ to imply that there exists a non-negative constant $c$ that depends on $a, b$ such that $A \leq cB, A = cB$ respectively. We write in a standard form the vorticity $w = \nabla \times v$ and current density $j = \nabla \times b$. 

Let us firstly obtain the conserved quantities. We have after taking $L^{2}$-inner products on (1a)-(1b) with $(u, b)$ respectively and integrating in time, 
\begin{align}
&\sup_{t \in [0,T]} (\lVert u(t) \rVert_{L^{2}}^{2} + \lVert \nabla u(t)\rVert_{L^{2}}^{2} + \lVert b(t)\rVert_{L^{2}}^{2})\nonumber\\
&+ \int_{0}^{T} \lVert \Lambda^{r_{1}} u\rVert_{L^{2}}^{2} + \lVert \Lambda^{r_{1}} \nabla u\rVert_{L^{2}}^{2} + \lVert \Lambda^{r_{2}} b\rVert_{L^{2}}^{2} \lesssim 1.
\end{align}
We now state some key lemmas: 

\begin{lemma}
(cf. [4]) Let $f$ be divergence-free vector field such that $\nabla f \in L^{p}, p \in (1,\infty)$. Then 
\begin{equation*}
\lVert \nabla f\rVert_{L^{p}} \leq c\frac{p^{2}}{p-1}\lVert \nabla\times  f\rVert_{L^{p}}.
\end{equation*}
\end{lemma}

\begin{lemma}
(cf. [17]) Let $f \in W^{\delta, p_{1}}(\mathbb{R}^{2}) \cap L^{q_{2}}(\mathbb{R}^{2}), g \in W^{\delta, p_{2}}(\mathbb{R}^{2}) \cap L^{q_{1}}(\mathbb{R}^{2}), \delta \geq 0, 1 < p_{k} < \infty, 1 < q_{k} \leq \infty, \frac{1}{p_{k}} + \frac{1}{q_{k}} = \frac{1}{p}, k = 1, 2$. Then 
\begin{equation*}
\lVert fg\rVert_{\dot{W}^{\delta, p}} \lesssim (\lVert f\rVert_{\dot{W}^{\delta, p_{1}}}\lVert g\rVert_{L^{q_{1}}} + \lVert f\rVert_{L^{q_{2}}} \lVert g\rVert_{\dot{W}^{\delta, p_{2}}}).
\end{equation*}
\end{lemma}

\begin{lemma}
(cf. [18]) Let $f,g$ be smooth such that $\nabla f \in L^{p_{1}}(\mathbb{R}^{2}), \Lambda^{s-1}g \in L^{p_{2}}(\mathbb{R}^{2}), \Lambda^{s}f \in L^{p_{3}}(\mathbb{R}^{2}), g \in L^{p_{4}}(\mathbb{R}^{2}), p \in (1,\infty), \frac{1}{p} = \frac{1}{p_{1}}+\frac{1}{p_{2}} = \frac{1}{p_{3}} + \frac{1}{p_{4}}, p_{2}, p_{3} \in (1, \infty), s > 0.$ Then 
\begin{equation*}
\lVert \Lambda^{s}(fg) - f\Lambda^{s}g\rVert_{L^{p}} \lesssim (\lVert \nabla f\rVert_{L^{p_{1}}}\lVert \Lambda^{s-1}g\rVert_{L^{p_{2}}} + \lVert \Lambda^{s}f\rVert_{L^{p_{3}}}\lVert g\rVert_{L^{p_{4}}}).
\end{equation*}
\end{lemma}

\begin{lemma}
(cf. [8], [16]) Let $\gamma \in [0, 1], x \in \mathbb{R}^{2}, \mathbb{T}^{2}, $ and $f, \Lambda^{2\gamma}f \in L^{p}, p \geq 2$. Then  
\begin{equation*}
2\int \lvert\Lambda^{\gamma}\lvert f\rvert^{\frac{p}{2}}\rvert^{2} \leq p \int\lvert f\rvert^{p-2}f\Lambda^{2\gamma}f.
\end{equation*}
\end{lemma}

\section{A priori estimates: case $\frac{1}{2} < r_{1} < 1, 0 < r_{2} < \frac{1}{2}$}

We fix $r_{1} \in \left(\frac{1}{2}, 1\right), r_{2} \in \left(0, \frac{1}{2}\right)$ and first obtain the following key proposition:

\begin{proposition}
Suppose $\nu, \eta > 0, \frac{1}{2} < r_{1} <1, 0 < r_{2} < \frac{1}{2}$ so that $r_{1} + r_{2} = 1, $ and $(v, b)$ solves (1a)-(1c) in $[0,T]$. Then $\forall \hspace{1mm} \gamma \in (1- r_{2}, 1)$, 
\begin{align}
&\sup_{t \in [0,T]} \left(\lVert \Lambda^{r_{1}}u(t)\rVert_{L^{2}}^{2} + \lVert \Lambda^{r_{1}}\nabla u(t)\rVert_{L^{2}}^{2} + \lVert \Lambda^{\gamma}b(t)\rVert_{L^{2}}^{2}\right)\\
&+ \int_{0}^{T} \lVert \Lambda^{2r_{1}}u\rVert_{L^{2}}^{2} + \lVert \Lambda^{2r_{1}}\nabla u\rVert_{L^{2}}^{2} + \lVert \Lambda^{r_{2} + \gamma}b\rVert_{L^{2}}^{2} d\tau \lesssim 1.\nonumber
\end{align}
\end{proposition}

\begin{proof}
We take $L^{2}$-inner products on (1a) with $\Lambda^{2r_{1}}u$ to obtain 
\begin{align}
&\frac{1}{2}\partial_{t} (\lVert \Lambda^{r_{1}} u\rVert_{L^{2}}^{2} + \lVert \Lambda^{r_{1}}\nabla u\rVert_{L^{2}}^{2}) + \lVert \Lambda^{2r_{1}} u\rVert_{L^{2}}^{2} + \lVert \Lambda^{2r_{1}} \nabla u\rVert_{L^{2}}^{2}\\
=& -\int (u\cdot\nabla) u\cdot\Lambda^{2r_{1}} u + \int (u\cdot\nabla) \Delta u \cdot\Lambda^{2r_{1}} u\nonumber\\
&- \sum_{k=1}^{2} \int u_{k}\nabla u_{k} \cdot \Lambda^{2r_{1}}u + \sum_{k=1}^{2} \int \Delta u_{k}\nabla u_{k} \cdot \Lambda^{2r_{1}}u\nonumber\\
&- \int \nabla \left(\pi + \frac{1}{2}\lvert b\rvert^{2}\right) \cdot\Lambda^{2r_{1}} u + \int (b\cdot\nabla) b \cdot \Lambda^{2r_{1}}u : = \sum_{i=1}^{6} I_{i}. \nonumber
\end{align}
On the other hand, for fixed $r_{2} > 0$ we find 
\begin{equation}
\gamma \in (1-r_{2}, 1)
\end{equation}
and take $L^{2}$-inner products of (1b) with $\Lambda^{2\gamma} b$ to obtain 
\begin{align}
&\frac{1}{2}\partial_{t} \lVert \Lambda^{\gamma}b\rVert_{L^{2}}^{2} + \lVert \Lambda^{r_{2} + \gamma}b\rVert_{L^{2}}^{2}\\
=& -\int \Lambda^{\gamma}[(u\cdot\nabla) b - u\cdot\nabla\Lambda^{\gamma}b]\cdot\Lambda^{\gamma} b\nonumber\\
& \hspace{12mm} + \int \Lambda^{\gamma - r_{2}}((b\cdot\nabla) u) \cdot\Lambda^{\gamma + r_{2} }b := II_{1} + II_{2}.\nonumber
\end{align}
We first estimate 
\begin{align}
I_{1} =& \int (u\cdot\nabla) \Lambda^{2r_{1}}u\cdot u
\leq \lVert u\rVert_{L^{2}} \lVert \Lambda^{2r_{1}}\nabla u\rVert_{L^{2}} \lVert u\rVert_{L^{\infty}}\\
\lesssim& \lVert \Lambda^{2r_{1}} \nabla u\rVert_{L^{2}} \lVert u\rVert_{L^{2}}^{1+\frac{r_{1}}{1+r_{1}}}\lVert \Lambda^{1+r_{1}} u\rVert_{L^{2}}^{\frac{1}{1+r_{1}}} \leq \frac{1}{8} \lVert \Lambda^{2r_{1}} \nabla u\rVert_{L^{2}}^{2} + c(1+ \lVert \Lambda^{r_{1}}\nabla u\rVert_{L^{2}}^{2})\nonumber
\end{align}
by integration by parts using the incompressibility of $u$, H$\ddot{o}$lder's, Gagliardo-Nirenberg and Young's inequalities and (2). Next, for fixed $r_{1} \in \left(\frac{1}{2}, 1\right)$, we find 
\begin{equation}
\epsilon \in \left(0, r_{1} - \frac{1}{2}\right)
\end{equation}
so that on $I_{2}$ from (4) we can estimate 
\begin{align}
I_{2} = -\int (u\cdot\nabla) \Lambda^{2r_{1}} u \cdot \Delta u
\leq \lVert u\rVert_{L^{\frac{1}{r_{1} - \frac{1}{2} - \epsilon}}}\lVert \Lambda^{2r_{1}} \nabla u\rVert_{L^{2}} \lVert \Delta u\rVert_{L^{\frac{1}{1-r_{1} + \epsilon}}}
\end{align}
by integration by parts and H$\ddot{o}$lder's inequalities. Now 
\begin{equation}
\lVert u\rVert_{L^{\frac{1}{r_{1} - \frac{1}{2} - \epsilon}}} \lesssim \lVert u\rVert_{L^{2}}^{2r_{1} - 1 - 2\epsilon}\lVert \nabla u\rVert_{L^{2}}^{2+2\epsilon - 2r_{1}} \lesssim 1
\end{equation}
by Gagliardo-Nirenberg inequality and (2). We also estimate similarly
\begin{equation}
\lVert \Delta u\rVert_{L^{\frac{1}{1 - r_{1} + \epsilon}}} \lesssim \lVert \nabla u\rVert_{L^{2}}^{\frac{\epsilon}{r_{1}}}\lVert \Lambda^{2r_{1}} \nabla u\rVert_{L^{2}}^{\frac{r_{1} -\epsilon}{r_{1}}} \lesssim \lVert \Lambda^{2r_{1}} \nabla u\rVert_{L^{2}}^{\frac{r_{1} -\epsilon}{r_{1}}} 
\end{equation}
by Gagliardo-Nirenberg inequality and (2). Taking into account of (10) and (11) in (9), we obtain 
\begin{equation}
I_{2} \lesssim \lVert \Lambda^{2r_{1}} \nabla u\rVert_{L^{2}}^{1+\frac{r_{1} - \epsilon}{r_{1}}} \leq \frac{1}{8} \lVert \Lambda^{2r_{1}}\nabla u\rVert_{L^{2}}^{2} + c
\end{equation}
by Young's inequality. Next, clearly $I_{3} = 0$ by incompressibility of $u$. Next, 
\begin{align}
I_{4} 
\lesssim \lVert \Delta u\rVert_{L^{\frac{1}{1-r_{1}}}} \lVert \nabla u\rVert_{L^{2}} \lVert \Lambda^{2r_{1}} u\rVert_{L^{\frac{2}{2r_{1} -1}}}
\end{align}
by H$\ddot{o}$lder's inequalities. We use Gagliardo-Nirenberg inequality and (2) to bound 
\begin{equation}
\lVert \Lambda^{2r_{1}} u\rVert_{L^{\frac{2}{2r_{1} - 1}}}  \lesssim \lVert \Lambda u\rVert_{L^{2}}^{\frac{2r_{1} - 1}{2r_{1}}} \lVert \Lambda^{2r_{1}} \Lambda u\rVert_{L^{2}}^{\frac{1}{2r_{1}}} \lesssim \lVert \Lambda^{2r_{1}} \nabla u\rVert_{L^{2}}^{\frac{1}{2r_{1}}}.
\end{equation}
Considering (14) and (2) in (13) and using the Sobolev embedding of $\dot{H}^{2r_{1} -1}(\mathbb{R}^{2})\hookrightarrow L^{\frac{1}{1-r_{1}}}(\mathbb{R}^{2})$ we obtain 
\begin{equation}
I_{4} \lesssim \lVert \Lambda^{2r_{1}} \nabla u\rVert_{L^{2}}^{1+ \frac{1}{2r_{1}}} \leq \frac{1}{8} \lVert \Lambda^{2r_{1}} \nabla u\rVert_{L^{2}}^{2} + c
\end{equation}
by Young's inequality. Due to the incompressibility of $u$, $I_{5} = 0$. Lastly, 
\begin{align}
I_{6} =& -\int (b\cdot\nabla) \Lambda^{2r_{1}} u \cdot b \leq \lVert b\rVert_{L^{\frac{2}{r_{2}}}} \lVert \nabla \Lambda^{2r_{1}} u\rVert_{L^{2}} \lVert b\rVert_{L^{\frac{2}{1-r_{2}}}}\\
\lesssim& \lVert \Lambda^{1-r_{2}} b\rVert_{L^{2}} \lVert \nabla \Lambda^{2r_{1}} u\rVert_{L^{2}} \lVert \Lambda^{r_{2}} b\rVert_{L^{2}}\nonumber\\
\leq& \frac{1}{8} \lVert \nabla \Lambda^{2r_{1}} u\rVert_{L^{2}}^{2} + c\lVert b\rVert_{L^{2}}^{2(\frac{\gamma + r_{2} - 1}{\gamma})}\lVert \Lambda^{\gamma} b\rVert_{L^{2}}^{2(\frac{1-r_{2}}{\gamma})}\lVert \Lambda^{r_{2}} b\rVert_{L^{2}}^{2}\nonumber\\
\leq& \frac{1}{8} \lVert \nabla \Lambda^{2r_{1}} u\rVert_{L^{2}}^{2} + c(1+ \lVert \Lambda^{\gamma} b\rVert_{L^{2}}^{2}) \lVert \Lambda^{r_{2}} b\rVert_{L^{2}}^{2}\nonumber
\end{align}
by integration by parts using incompressibility of $b$, H$\ddot{o}$lder's inequality, the Sobolev embeddings of $\dot{H}^{1-r_{2}}(\mathbb{R}^{2}) \hookrightarrow L^{\frac{2}{r_{2}}} (\mathbb{R}^{2})$ and $\dot{H}^{r_{2}} (\mathbb{R}^{2})\hookrightarrow L^{\frac{2}{1-r_{2}}} (\mathbb{R}^{2})$, Young's and Gagliardo-Nirenberg inequalities and (2). From (4), (7), (12), (15), (16) we thus have after absorbing the dissipative terms 
\begin{align}
&\partial_{t} (\lVert \Lambda^{r_{1}} u\rVert_{L^{2}}^{2} + \lVert \Lambda^{r_{1}} \nabla u\rVert_{L^{2}}^{2}) + \lVert \Lambda^{2r_{1}} u\rVert_{L^{2}}^{2} + \lVert \Lambda^{2r_{1}}\nabla u\rVert_{L^{2}}^{2}\\
\lesssim& (1+ \lVert \Lambda^{r_{1}}\nabla u\rVert_{L^{2}}^{2}) + (1+ \lVert \Lambda^{\gamma} b\rVert_{L^{2}}^{2}) \lVert \Lambda^{r_{2}} b\rVert_{L^{2}}^{2}.\nonumber
\end{align}

Next, we work on (6): firstly, 
\begin{align}
II_{1} \lesssim& \left(\lVert \nabla u\rVert_{L^{\frac{2}{r_{2}}}} \lVert \Lambda^{\gamma} b\rVert_{L^{\frac{2}{1-r_{2}}}} + \lVert \Lambda^{\gamma} u\rVert_{L^{\frac{2}{\gamma + r_{2} - 1}}}\lVert \nabla b\rVert_{L^{\frac{2}{2-\gamma - r_{2}}}}\right) \lVert \Lambda^{\gamma} b\rVert_{L^{2}}\\
\lesssim& (\lVert \Lambda^{2-r_{2}} u\rVert_{L^{2}} \lVert \Lambda^{\gamma + r_{2}} b\rVert_{L^{2}} + \lVert \Lambda^{2-r_{2}} u\rVert_{L^{2}} \lVert \Lambda^{\gamma + r_{2}}b\rVert_{L^{2}}) \lVert \Lambda^{\gamma} b\rVert_{L^{2}}\nonumber\\
\leq& \frac{1}{4} \lVert \Lambda^{\gamma + r_{2}} b\rVert_{L^{2}}^{2} + c\lVert \Lambda^{r_{1}} \nabla u\rVert_{L^{2}}^{2} \lVert \Lambda^{\gamma} b\rVert_{L^{2}}^{2}\nonumber
\end{align}
by H$\ddot{o}$lder's inequality and Lemma 2.3, the Sobolev embeddings of $\dot{H}^{1-r_{2}}(\mathbb{R}^{2})\hookrightarrow L^{\frac{2}{r_{2}}}(\mathbb{R}^{2}), \dot{H}^{r_{2}}(\mathbb{R}^{2}) \hookrightarrow L^{\frac{2}{1-r_{2}}}(\mathbb{R}^{2}), \dot{H}^{2-\gamma-r_{2}}(\mathbb{R}^{2})\hookrightarrow L^{\frac{2}{\gamma + r_{2} - 1}}(\mathbb{R}^{2})$ and $\dot{H}^{\gamma + r_{2} -1}(\mathbb{R}^{2}) \hookrightarrow L^{\frac{2}{2-\gamma-r_{2}}}(\mathbb{R}^{2})$; we also used Young's inequality. 

Finally, we work on $II_{2}$ from (6):
\begin{align}
II_{2} \lesssim& (\lVert b\rVert_{L^{\frac{2}{1-\gamma}}}\lVert \Lambda^{\gamma - r_{2}}\nabla u\rVert_{L^{\frac{2}{\gamma}}} + \lVert \Lambda^{\gamma - r_{2}}b\rVert_{L^{\frac{2}{r_{1}}}}\lVert \nabla u\rVert_{L^{\frac{2}{1-r_{1}}}})\lVert \Lambda^{\gamma + r_{2}} b\rVert_{L^{2}}\\
\lesssim& (\lVert \Lambda^{\gamma} b\rVert_{L^{2}} \lVert \Lambda^{r_{1}} \nabla u\rVert_{L^{2}} + \lVert \Lambda^{\gamma}b\rVert_{L^{2}} \lVert \Lambda^{r_{1}} \nabla u\rVert_{L^{2}}) \lVert \Lambda^{\gamma + r_{2}} b\rVert_{L^{2}}\nonumber\\
\leq& \frac{1}{4} \lVert \Lambda^{\gamma + r_{2}} b\rVert_{L^{2}}^{2} + c\lVert \Lambda^{\gamma} b\rVert_{L^{2}}^{2} \lVert \Lambda^{r_{1}} \nabla u\rVert_{L^{2}}^{2}\nonumber
\end{align}
by H$\ddot{o}$lder's inequality, Lemma 2.2, the Sobolev embeddings of $\dot{H}^{\gamma}(\mathbb{R}^{2}) \hookrightarrow L^{\frac{2}{1-\gamma}}(\mathbb{R}^{2})$, $\dot{H}^{1-\gamma}(\mathbb{R}^{2}) \hookrightarrow L^{\frac{2}{\gamma}}(\mathbb{R}^{2}), \dot{H}^{1-r_{1}}(\mathbb{R}^{2})\hookrightarrow L^{\frac{2}{r_{1}}}(\mathbb{R}^{2})$ and $\dot{H}^{r_{1}}(\mathbb{R}^{2})\hookrightarrow L^{\frac{2}{1-r_{1}}}(\mathbb{R}^{2})$, and Young's inequality. 

From (6), (18) and (19), absorbing the diffusive term we obtain 
\begin{equation}
\partial_{t} \lVert \Lambda^{\gamma} b\rVert_{L^{2}}^{2} + \lVert \Lambda^{\gamma + r_{2}} b\rVert_{L^{2}}^{2}
\lesssim \lVert \Lambda^{r_{1}}\nabla u\rVert_{L^{2}}^{2} \lVert \Lambda^{\gamma} b\rVert_{L^{2}}^{2}.
\end{equation}
Summing (17) and (20) and using (2), Gronwall's inequality completes the proof of Proposition 3.1. 
\end{proof}

We make use of the vorticity formulation and obtain higher regularity:

\begin{proposition}
Suppose $\nu, \eta > 0, \frac{1}{2} < r_{1} < 1, 0 < r_{2} < \frac{1}{2}$ so that $r_{1} + r_{2} = 1, $ and $(v, b)$ solves (1a)-(1c) in $[0,T]$. Then with $w = \nabla \times v$
\begin{align}
\sup_{t \in [0,T]} \left(\lVert w(t)\rVert_{L^{2}}^{2} + \lVert \Lambda^{1+r_{1}} b(t)\rVert_{L^{2}}^{2}\right) + \int_{0}^{T} \lVert \Lambda^{r_{1}}w\rVert_{L^{2}}^{2} + \lVert \Lambda^{2}b\rVert_{L^{2}}^{2}  d\tau \lesssim 1.\nonumber
\end{align}
\end{proposition}

\begin{proof}
We apply $\nabla \times $ on (1a) and obtain 
\begin{equation}
\partial_{t} w + (u\cdot\nabla) w + \Lambda^{2r_{1}} w = (b\cdot\nabla) j. 
\end{equation}
We take $L^{2}$-inner products on (21) with $w$ to obtain 
\begin{equation}
\frac{1}{2}\partial_{t} \lVert w\rVert_{L^{2}}^{2} + \lVert \Lambda^{r_{1}} w\rVert_{L^{2}}^{2} = \int (b\cdot\nabla) j w \leq \lVert \nabla j\rVert_{L^{2}} \lVert b\rVert_{L^{\infty}} \lVert w\rVert_{L^{2}} 
\end{equation}
by H$\ddot{o}$lder's inequality. Now for fixed $r_{2} \in \left(0, \frac{1}{2}\right)$, we can find 
\begin{equation}
\delta \in (1, 1+r_{2})
\end{equation}
so that by Proposition 3.1, 
\begin{equation}
\int_{0}^{T} \lVert \Lambda^{\delta} b\rVert_{L^{2}}^{2} d\tau \lesssim 1.
\end{equation}
We now use Gagliardo-Nirenberg and Young's inequalities and (2) to estimate the right hand side of (22) by 
\begin{equation}
c\lVert \Lambda^{2} b\rVert_{L^{2}}\lVert b\rVert_{L^{2}}^{\frac{\delta - 1}{\delta}}\lVert \Lambda^{\delta} b\rVert_{L^{2}}^{\frac{1}{\delta}} \lVert w\rVert_{L^{2}}\leq \frac{1}{4} \lVert \Lambda^{2} b\rVert_{L^{2}}^{2} + c\left(1 + \lVert \Lambda^{\delta} b\rVert_{L^{2}}^{2}\right) \lVert w\rVert_{L^{2}}^{2}.
\end{equation}
On the other hand, taking $L^{2}$-inner products on (1b) with $\Lambda^{2+2r_{1}}b$, we obtain 
\begin{align}
\frac{1}{2}\partial_{t} \lVert \Lambda^{1+r_{1}}b\rVert_{L^{2}}^{2} + \lVert \Lambda^{2} b\rVert_{L^{2}}^{2} =& -\int (u\cdot\nabla) b\cdot\Lambda^{2+2r_{1}} b + \int (b\cdot\nabla) u \cdot \Lambda^{2+2r_{1}} b\\
:=& III_{1} + III_{2}\nonumber 
\end{align}
where we used that $r_{1} + r_{2} = 1$. We rewrite $III_{1}$ using the incompressibility of $u$ and estimate 
\begin{align}
III_{1} =& \int [\Lambda^{1+r_{1}}((u\cdot\nabla) b) - u\cdot\nabla \Lambda^{1+r_{1}} b]\cdot\Lambda^{1+r_{1}} b\\
\lesssim& (\lVert \nabla u\rVert_{L^{\frac{2}{r_{2}}}}\lVert \Lambda^{1+r_{1}} b\rVert_{L^{\frac{2}{1-r_{2}}}} + \lVert \Lambda^{1+r_{1}} u\rVert_{L^{\frac{2}{\delta -1}}} \lVert \nabla b\rVert_{L^{\frac{2}{2-\delta}}})\Vert \Lambda^{1+r_{1}} b\rVert_{L^{2}}\nonumber\\
\lesssim& \left(\lVert \Lambda^{1+r_{1}} u\rVert_{L^{2}} \lVert \Lambda^{2} b\rVert_{L^{2}}  + \lVert \Lambda^{1+r_{1}} u\rVert_{L^{2}}^{\frac{\delta}{2}} \lVert \Lambda^{3+r_{1}} u\rVert_{L^{2}}^{\frac{2-\delta}{2}}\lVert \Lambda^{\delta} b\rVert_{L^{2}}\right) \lVert \Lambda^{1+r_{1}} b\rVert_{L^{2}}\nonumber
\end{align}
where we used H$\ddot{o}$lder's inequality, Lemma 2.3, the Sobolev embeddings of $\dot{H}^{1-r_{2}}(\mathbb{R}^{2}) \hookrightarrow L^{\frac{2}{r_{2}}} (\mathbb{R}^{2}), \dot{H}^{r_{2}}(\mathbb{R}^{2}) \hookrightarrow L^{\frac{2}{1-r_{2}}}(\mathbb{R}^{2})$ and $\dot{H}^{\delta -1}(\mathbb{R}^{2}) \hookrightarrow L^{\frac{2}{2-\delta}}(\mathbb{R}^{2})$, and Gagliardo-Nirenberg inequality. We bound (27) using Proposition 3.1 and Young's inequalities: 
\begin{align}
III_{1} \lesssim&  \lVert \Lambda^{2} b\rVert_{L^{2}} \lVert \Lambda^{1+r_{1}} b\rVert_{L^{2}} +  \lVert \Lambda^{3+r_{1}} u\rVert_{L^{2}}^{\frac{2-\delta}{2}}\lVert \Lambda^{\delta} b\rVert_{L^{2}} \lVert \Lambda^{1+r_{1}} b\rVert_{L^{2}}\\
\leq& \frac{1}{4} \left( \lVert \Lambda^{2} b\rVert_{L^{2}}^{2} +  \lVert \Lambda^{r_{1}} w\rVert_{L^{2}}^{2}\right) + c(1+ \lVert \Lambda^{1+r_{1}} b\rVert_{L^{2}}^{2})(1+ \lVert \Lambda^{\delta} b\rVert_{L^{2}}^{2}).\nonumber
\end{align}

Next, we employ H$\ddot{o}$lder's inequality, Lemma 2.2, Gagliardo-Nirenberg inequality with same $\delta$ defined in (23), (2) and Young's inequalities to obtain 
\begin{align}
III_{2} =& \int \Lambda^{1+r_{1}} ((b\cdot\nabla) u) \cdot \Lambda^{1+r_{1}} b\\
\lesssim& (\lVert b\rVert_{L^{\infty}} \lVert \Lambda^{1+r_{1}}\nabla u\rVert_{L^{2}} + \lVert \Lambda^{1+r_{1}}b\rVert_{L^{2}} \lVert \nabla u\rVert_{L^{\infty}})\lVert \Lambda^{1+r_{1}} b\rVert_{L^{2}}\nonumber\\
\lesssim& (\lVert b\rVert_{L^{2}}^{\frac{\delta - 1}{\delta}}\lVert \Lambda^{\delta} b\rVert_{L^{2}}^{\frac{1}{\delta}}\lVert \Lambda^{2+r_{1}} u\rVert_{L^{2}}\nonumber\\
& \hspace{5mm} + \lVert \Lambda^{1+r_{1}} b\rVert_{L^{2}}\lVert \nabla u\rVert_{L^{2}}^{\frac{2r_{1} - 1}{2r_{1}}}\lVert \Lambda^{2r_{1}}\nabla u\rVert_{L^{2}}^{\frac{1}{2r_{1}}}) \lVert \Lambda^{1+r_{1}} b\rVert_{L^{2}}\nonumber\\
\leq& \frac{1}{4} \lVert \Lambda^{r_{1}} w\rVert_{L^{2}}^{2} + c(1+ \lVert \Lambda^{\delta} b\rVert_{L^{2}}^{\frac{2}{\delta}} + \lVert \Lambda^{2r_{1}} \nabla u\rVert_{L^{2}}^{\frac{1}{2r_{1}}})(1+ \lVert \Lambda^{1+r_{1}} b\rVert_{L^{2}}^{2})\nonumber\\
\leq& \frac{1}{4} \lVert \Lambda^{r_{1}} w\rVert_{L^{2}}^{2} + c(1+ \lVert \Lambda^{\delta} b\rVert_{L^{2}}^{2} + \lVert \Lambda^{2r_{1}} \nabla u\rVert_{L^{2}}^{2})(1+ \lVert \Lambda^{1+r_{1}} b\rVert_{L^{2}}^{2})\nonumber
\end{align}
where we used that $\lVert \Lambda^{2+r_{1}}u\rVert_{L^{2}} \lesssim  \lVert \Lambda^{r_{1}} w\rVert_{L^{2}}^{2}$. Considering (22), (25), (26), (28), (29) we obtain after absorbing the dissipative and diffusive terms
\begin{align*}
&\partial_{t} \left(\lVert w\rVert_{L^{2}}^{2} + \lVert \Lambda^{1+r_{1}} b\rVert_{L^{2}}^{2}\right) + \lVert \Lambda^{r_{1}} w\rVert_{L^{2}}^{2} + \lVert \Lambda^{2} b\rVert_{L^{2}}^{2} \\
\lesssim& \left(1 + \lVert \Lambda^{\delta} b\rVert_{L^{2}}^{2} + \lVert \Lambda^{2r_{1}}\nabla u\rVert_{L^{2}}^{2}\right) (1 + \lVert w\rVert_{L^{2}}^{2} + \lVert \Lambda^{1+r_{1}} b\rVert_{L^{2}}^{2}).
\end{align*}
Due to Proposition 3.1 and (24), Gronwall's inequality completes the proof of Proposition 3.2. 
\end{proof}

\begin{proposition}
Suppose $\nu, \eta > 0, \frac{1}{2} < r_{1} < 1, 0 < r_{2} < \frac{1}{2}$ so that $r_{1} + r_{2} = 1, $ and $(v, b)$ solves (1a)-(1c) in $[0,T]$. Then 
\begin{align}
\sup_{t \in [0,T]} \lVert \Lambda^{3}b(t)\rVert_{L^{2}}^{2}  + \int_{0}^{T} \lVert \Lambda^{3 + r_{2}}b\rVert_{L^{2}}^{2} d\tau \lesssim 1.\nonumber
\end{align}
\end{proposition}

\begin{proof}
We take $L^{2}$-inner products of (1b) with $\Lambda^{6} b$ to obtain 
\begin{align}
&\frac{1}{2}\partial_{t}\lVert \Lambda^{3} b\rVert_{L^{2}}^{2} + \lVert \Lambda^{3+r_{2}} b\rVert_{L^{2}}^{2}\\
=& -\int [\Lambda^{3}((u\cdot\nabla) b) - u\cdot\nabla\Lambda^{3} b] \cdot \Lambda^{3} b\nonumber\\
& \hspace{10mm}  + \int \Lambda^{3-r_{2}} ((b\cdot\nabla)u)\cdot\Lambda^{3+r_{2}} b := IV_{1} + IV_{2}.\nonumber
\end{align}

We use H$\ddot{o}$lder's inequality, Lemma 2.3, the Sobolev embedding of $H^{2}(\mathbb{R}^{2})\hookrightarrow L^{\infty}(\mathbb{R}^{2})$, $\dot{H}^{r_{1}}(\mathbb{R}^{2}) \hookrightarrow L^{\frac{2}{1-r_{1}}}(\mathbb{R}^{2})$ and $\dot{H}^{1-r_{1}}(\mathbb{R}^{2}) \hookrightarrow L^{\frac{2}{r_{1}}}(\mathbb{R}^{2})$, (2) and Gagliardo-Nirenberg inequality to obtain 
\begin{align*}
IV_{1} \lesssim& (\lVert \nabla u\rVert_{L^{\infty}} \lVert \Lambda^{3} b\rVert_{L^{2}} + \lVert \Lambda^{3} u\rVert_{L^{\frac{2}{r_{1}}}}\lVert \nabla b\rVert_{L^{\frac{2}{1-r_{1}}}} ) \lVert \Lambda^{3} b\rVert_{L^{2}}\\
\lesssim& (\lVert \nabla u\rVert_{L^{2}} + \lVert \Lambda^{3} u\rVert_{L^{2}})\lVert \Lambda^{3} b\rVert_{L^{2}}^{2} + \lVert \Lambda^{4-r_{1}}u\rVert_{L^{2}} \lVert \Lambda^{1+r_{1}} b\rVert_{L^{2}}\lVert \Lambda^{3} b\rVert_{L^{2}}\\
\lesssim& (1+ \lVert w\rVert_{L^{2}}^{2})\lVert \Lambda^{3} b\rVert_{L^{2}}^{2}+ \lVert u\rVert_{L^{2}}^{\frac{2r_{1} - 1}{3+r_{1}}}\lVert \Lambda^{3+r_{1}}u\rVert_{L^{2}}^{\frac{4-r_{1}}{3+r_{1}}}\lVert \Lambda^{1+r_{1}} b\rVert_{L^{2}} \lVert \Lambda^{3} b\rVert_{L^{2}}.
\end{align*}
By (2) and Proposition 3.2, this gives us the bound of 
\begin{equation}
IV_{1} \lesssim (1+ \lVert \Lambda^{3} b\rVert_{L^{2}}^{2})(1+ \lVert \Lambda^{3+r_{1}} u\rVert_{L^{2}}^{2}).
\end{equation}

Next, we use H$\ddot{o}$lder's inequality, Lemma 2.2, Gagliardo-Nirenberg inequalities, the Sobolev embedding of $H^{2}(\mathbb{R}^{2}) \hookrightarrow L^{\infty}(\mathbb{R}^{2})$, (2), Proposition 3.2, and Young's inequalities to obtain
\begin{align}
IV_{2} \lesssim& (\lVert \Lambda^{3-r_{2}} b\rVert_{L^{2}} \lVert \nabla u\rVert_{L^{\infty}} + \lVert b\rVert_{L^{\infty}}\lVert \Lambda^{3-r_{2}} \nabla u\rVert_{L^{2}}) \lVert \Lambda^{3+r_{2}} b\rVert_{L^{2}}\\
\lesssim& (\lVert \Lambda^{1+r_{1}} b\rVert_{L^{2}}^{\frac{r_{2}}{2-r_{1}}}\lVert \Lambda^{3} b\rVert_{L^{2}}^{\frac{2-r_{1} - r_{2}}{2-r_{1}}} (\lVert \nabla u\rVert_{L^{2}}+ \lVert \Lambda^{3} u\rVert_{L^{2}})\nonumber\\
&+ \lVert b\rVert_{L^{2}}^{\frac{r_{1}}{1+r_{1}}}\lVert \Lambda^{1+r_{1}} b\rVert_{L^{2}}^{\frac{1}{1+r_{1}}} \lVert \Lambda^{3+r_{1}} u\rVert_{L^{2}}) \lVert \Lambda^{3+r_{2}} b\rVert_{L^{2}}\nonumber\\
\lesssim& (\lVert \Lambda^{3} b\rVert_{L^{2}}^{\frac{1}{2-r_{1}}}(1+ \lVert w\rVert_{L^{2}}) + \lVert \Lambda^{3+r_{1}} u\rVert_{L^{2}}) \lVert \Lambda^{3+r_{2}} b\rVert_{L^{2}}\nonumber\\
\leq& \frac{1}{2} \lVert \Lambda^{3+r_{2}}b\rVert_{L^{2}}^{2} + c(1+ \lVert \Lambda^{3+r_{1}}u\rVert_{L^{2}}^{2})(1+ \lVert \Lambda^{3} b\rVert_{L^{2}}^{2}).\nonumber
\end{align}
Considering (30), (31), (32), after absorbing the diffusive term, Gronwall's inequality using Proposition 3.2 completes the proof of Proposition 3.3. 
\end{proof}

\begin{proposition}
Suppose $\nu, \eta > 0, \frac{1}{2} < r_{1} < 1, 0 < r_{2} < \frac{1}{2}$ so that $r_{1} + r_{2} = 1, $ and $(v, b)$ solves (1a)-(1c) in $[0,T]$. Then for any $p \in [2, \infty]$, with $ w = \nabla \times v$
\begin{equation}
\sup_{t \in [0,T]} \lVert w(t)\rVert_{L^{p}} \lesssim 1.\nonumber
\end{equation}
\end{proposition}

\begin{proof}
For any $p \in (2, \infty)$, we multiply (21) by $\lvert w\rvert^{p-2}w$, integrate in space to obtain by H$\ddot{o}$lder's inequality 
\begin{equation*}
\frac{1}{p} \partial_{t} \lVert w\rVert_{L^{p}}^{p} + \int \Lambda^{2r_{1}} w\lvert w\rvert^{p-2} w  \leq \lVert b\rVert_{L^{2p}} \lVert \nabla j\rVert_{L^{2p}} \lVert w\rVert_{L^{p}}^{p-1}.
\end{equation*}
We use Lemma 2.4 to bound the dissipative term from below so that 
\begin{equation*}
\partial_{t} \lVert w\rVert_{L^{p}} \leq \lVert b\rVert_{L^{2p}} \lVert \nabla j\rVert_{L^{2p}}.
\end{equation*}
Taking $p \to \infty$, integrating in time, we obtain with $w(0) = \nabla \times v_{0}$, 
\begin{align*}
&\sup_{t \in [0,T]} \lVert w(t)\rVert_{L^{\infty}}\\
\leq& \lVert w(0)\rVert_{L^{\infty}} + \sup_{t \in [0,T]} \lVert b(t)\rVert_{L^{\infty}} \int_{0}^{T} \lVert \nabla j(\tau)\rVert_{L^{\infty}} d\tau\\
\lesssim& 1 + \sup_{t \in [0,T]} (\lVert b(t)\rVert_{L^{2}}+ \lVert \Lambda^{3} b(t)\rVert_{L^{2}})\lVert \Lambda^{r_{1}} j(t)\rVert_{L^{2}}^{\frac{r_{2}}{2+ r_{2} - r_{1}}}\int_{0}^{T} \lVert \Lambda^{2+r_{2}} j\rVert_{L^{2}}^{\frac{2-r_{1}}{2+ r_{2} - r_{1}}}d\tau
\end{align*}
by the Sobolev embedding of $H^{3}(\mathbb{R}^{2})  \hookrightarrow L^{\infty}(\mathbb{R}^{2})$ and Gagliardo-Nirenberg inequality. The bounds from (2), Propositions 3.2 and 3.3 imply that the right hand side is finite. Interpolating between $p \in [2, \infty]$ using Proposition 3.2 completes the proof of Proposition 3.4. 
\end{proof}

\begin{proposition}
Suppose $\nu, \eta > 0, \frac{1}{2} < r_{1} < 1, 0 < r_{2} < \frac{1}{2}$ so that $r_{1} + r_{2} = 1, $ and $(v, b)$ solves (1a)-(1c) in $[0,T]$. Then 
\begin{equation}
\sup_{t \in [0,T]} \lVert \Lambda^{3} v(t)\rVert_{L^{2}} + \int_{0}^{T} \lVert \Lambda^{3+r_{1}} v\rVert_{L^{2}}^{2} d\tau \lesssim 1.\nonumber
\end{equation}
\end{proposition}

\begin{proof}
We take $L^{2}$-inner products of (1a) with $\Lambda^{6} v$ and estimate 
\begin{align*}
&\frac{1}{2} \partial_{t} \lVert \Lambda^{3} v\rVert_{L^{2}}^{2} + \lVert \Lambda^{3+r_{1}}v\rVert_{L^{2}}^{2}\\
=& -\int [\Lambda^{3} ((u\cdot\nabla) v) - u\cdot\nabla \Lambda^{3} v] \cdot\Lambda^{3} v\\
& - \int \Lambda^{3} \left(\sum_{k=1}^{2} v_{k}\nabla u_{k}\right) \cdot\Lambda^{3} v + \int \Lambda^{3-r_{1}} \text{div} (b\otimes b) \cdot\Lambda^{3+r_{1}} v\\
\leq& \lVert \Lambda^{3} ((u\cdot\nabla) v) - u\cdot\nabla \Lambda^{3} v\rVert_{L^{\frac{2}{1+r_{1}}}} \lVert \Lambda^{3} v\rVert_{L^{\frac{2}{1-r_{1}}}}\\
&+ \left\lVert \Lambda^{3} \left(\sum_{k=1}^{2} v_{k}\nabla u_{k}\right)\right\rVert_{L^{\frac{2}{1+r_{1}}}} \lVert \Lambda^{3} v\rVert_{L^{\frac{2}{1-r_{1}}}} + \lVert \Lambda^{4-r_{1}} (b\otimes b)\rVert_{L^{2}} \lVert \Lambda^{3+r_{1}} v\rVert_{L^{2}}\\
\lesssim& (\lVert \nabla u\rVert_{L^{\frac{2}{r_{1}}}} \lVert \Lambda^{3} v\rVert_{L^{2}} + \lVert \Lambda^{3} u\rVert_{L^{2}}\lVert \nabla v\rVert_{L^{\frac{2}{r_{1}}}} ) \lVert \Lambda^{3 + r_{1}} v\rVert_{L^{2}}\\
&+ (\lVert \Lambda^{3} v\rVert_{L^{2}} \lVert \nabla u\rVert_{L^{\frac{2}{r_{1}}}} + \lVert v\rVert_{L^{\frac{2}{r_{1}}}} \lVert \Lambda^{3} \nabla u\rVert_{L^{2}}) \lVert \Lambda^{3+r_{1}} v\rVert_{L^{2}}\\
&+ \lVert b\rVert_{L^{\infty}} \lVert \Lambda^{3+r_{2}}b\rVert_{L^{2}} \lVert \Lambda^{3+r_{1}} v\rVert_{L^{2}}
\end{align*}
by H$\ddot{o}$lder's inequalities, Lemmas 2.3 and 2.2 and the Sobolev embedding of $\dot{H}^{r_{1}}(\mathbb{R}^{2}) \hookrightarrow L^{\frac{2}{1-r_{1}}}(\mathbb{R}^{2})$. We furthermore bound this by 
\begin{align*}
&\lVert u\rVert_{L^{2}}^{\frac{r_{1} + 1}{3}}\lVert \Lambda^{3} u\rVert_{L^{2}}^{\frac{2-r_{1}}{3}} \lVert \Lambda^{3} v\rVert_{L^{2}} \lVert \Lambda^{3+r_{1}} v\rVert_{L^{2}} + \lVert \Lambda^{3} v\rVert_{L^{2}}\lVert w\rVert_{L^{\frac{2}{r_{1}}}}\lVert \Lambda^{3+r_{1}} v\rVert_{L^{2}}\\
&+ (\lVert u\rVert_{L^{2}} + \lVert \Delta u\rVert_{L^{2}})^{r_{1}} \lVert w\rVert_{L^{2}}^{1-r_{1}}\lVert \Lambda^{3} v\rVert_{L^{2}} \lVert \Lambda^{3+r_{1}} v\rVert_{L^{2}}\\
&+ (\lVert b\rVert_{L^{2}} + \lVert \Lambda^{3} b\rVert_{L^{2}}) \lVert \Lambda^{3+r_{2}} b\rVert_{L^{2}} \lVert \Lambda^{3+r_{1}} v\rVert_{L^{2}}\\
\lesssim&  (\lVert w\rVert_{L^{2}} + 1)\lVert \Lambda^{3} v\rVert_{L^{2}} \lVert \Lambda^{3+r_{1}} v\rVert_{L^{2}} + \lVert \Lambda^{3+r_{2}} b\rVert_{L^{2}} \lVert \Lambda^{3+r_{1}} v\rVert_{L^{2}}\\
\leq& \frac{1}{2} \lVert \Lambda^{3+r_{1}} v\rVert_{L^{2}}^{2}+ c (\lVert \Lambda^{3} v\rVert_{L^{2}}^{2} + \lVert \Lambda^{3+r_{2}} b\rVert_{L^{2}}^{2})
\end{align*}
due to Gagliardo-Nirenberg inequality, Lemma 2.1 that implies 
\begin{equation*}
\lVert \nabla v\rVert_{L^{\frac{2}{r_{1}}}} \lesssim \lVert \nabla \times v\rVert_{L^{\frac{2}{r_{1}}}} \approx \lVert w\rVert_{L^{\frac{2}{r_{1}}}},
\end{equation*}
the Sobolev embedding of $H^{3}(\mathbb{R}^{2}) \hookrightarrow L^{\infty}(\mathbb{R}^{2})$, (2), Propositions 3.3 and 3.4. After absorbing the dissipative term, Gronwall's inequality using Proposition 3.3 completes the proof of Proposition 3.5. 
\end{proof}

\section{A priori estimates case $0 < r_{1} < \frac{1}{2}, \frac{1}{2} < r_{2} < 1$}

In the case $r_{2} > r_{1}$, because only the velocity is filtered, the relatively stronger diffusivity leads to a better balance. In fact, we can obtain the following proposition estimating only on $b$: 
\begin{proposition}
Suppose $\nu, \eta > 0, 0 < r_{1} < \frac{1}{2}, \frac{1}{2} <  r_{2} < 1$ so that $r_{1} + r_{2} = 1, $ and $(v, b)$ solves (1a)-(1c) in $[0,T]$. Then 
\begin{equation*}
\sup_{t \in [0,T]} \lVert \Lambda^{r_{2}} b(t)\rVert_{L^{2}} + \int_{0}^{T} \lVert \Lambda^{2r_{2}} b\rVert_{L^{2}}^{2} d\tau \lesssim 1.
\end{equation*}
\end{proposition}

\begin{proof}
We take $L^{2}$-inner products on (1b) with $\Lambda^{2r_{2}}b$ to estimate 
\begin{align}
\frac{1}{2} \partial_{t} \lVert \Lambda^{r_{2}} b\rVert_{L^{2}}^{2} + \lVert \Lambda^{2r_{2}} b\rVert_{L^{2}}^{2}
=& -\int (u\cdot\nabla) b \cdot \Lambda^{2 r_{2} } b - (b\cdot\nabla) u \cdot \Lambda^{2r_{2}} b\\
\leq& (\lVert \text{div} (u\otimes b) \rVert_{L^{2}} + \lVert \text{div} (b\otimes u) \rVert_{L^{2}}) \lVert \Lambda^{2r_{2}} b\rVert_{L^{2}}\nonumber
\end{align}
by H$\ddot{o}$lder's inequalities. Now for fixed $r_{2} \in \left(\frac{1}{2}, 1\right)$, we find 
\begin{equation}
\epsilon \in \left(0, \frac{2r_{2} -1}{2}\right)
\end{equation}
so that we denote for clarity 
\begin{equation*}
\lambda_{1} = \frac{1}{(\frac{2r_{2} - 1}{2}- \epsilon)}, \hspace{3mm}  \lambda_{2} = \frac{1}{(1-r_{2}) + \epsilon}
\end{equation*}
and bound (33) by Lemma 2.2, Gagliardo-Nirenberg inequalities, the Sobolev embeddings of $\dot{H}^{r_{1}}(\mathbb{R}^{2}) \hookrightarrow L^{\frac{2}{1-r_{1}}}(\mathbb{R}^{2})$ and $\dot{H}^{1-r_{1}}(\mathbb{R}^{2})\hookrightarrow L^{\frac{2}{r_{1}}}(\mathbb{R}^{2})$, (2) and Young's inequality as follows: 
\begin{align}
&\frac{1}{2} \partial_{t} \lVert \Lambda^{r_{2}} b\rVert_{L^{2}}^{2} + \lVert \Lambda^{2r_{2}} b\rVert_{L^{2}}^{2}\\
\lesssim& (\lVert u\rVert_{L^{\lambda_{1}}}\lVert \Lambda b\rVert_{L^{\lambda_{2}}} + \lVert \Lambda u\rVert_{L^{\frac{2}{1-r_{1}}}}\lVert b\rVert_{L^{\frac{2}{r_{1}}}} ) \lVert \Lambda^{2r_{2}} b\rVert_{L^{2}}\nonumber\\
\lesssim& (\lVert u\rVert_{L^{2}}^{2r_{2} - 1 - 2\epsilon}\lVert \Lambda u\rVert_{L^{2}}^{1-(2r_{2} - 1 - 2\epsilon)}\lVert b\rVert_{L^{2}}^{\frac{\epsilon}{r_{2}}}\lVert \Lambda^{2r_{2}}b\rVert_{L^{2}}^{1-\frac{\epsilon}{r_{2}}}\nonumber\\
& \hspace{40mm} + \lVert \Lambda^{1+r_{1}} u\rVert_{L^{2}} \lVert \Lambda^{1-r_{1}} b\rVert_{L^{2}}) \lVert \Lambda^{2r_{2}} b\rVert_{L^{2}}\nonumber\\
\lesssim& (\lVert \Lambda^{2r_{2}}b\rVert_{L^{2}}^{1-\frac{\epsilon}{r_{2}}} + \lVert \Lambda^{1+r_{1}} u\rVert_{L^{2}} \lVert \Lambda^{r_{2}} b\rVert_{L^{2}})\lVert \Lambda^{2r_{2}} b\rVert_{L^{2}}\nonumber\\
\leq& \frac{1}{2} \lVert \Lambda^{2r_{2}} b\rVert_{L^{2}}^{2} + c(1 + \lVert \Lambda^{r_{2}}b\rVert_{L^{2}}^{2})(1+ \lVert \Lambda^{r_{1}}\nabla u\rVert_{L^{2}}^{2}).\nonumber
\end{align}
After absorbing the diffusive term, Gronwall's inequality with (2) completes the proof of Proposition 4.1. 
\end{proof}

\begin{proposition}
Suppose $\nu, \eta > 0, 0 < r_{1} < \frac{1}{2}, \frac{1}{2} < r_{2} > 1$ so that $r_{1} + r_{2} = 1, $ and $(v, b)$ solves (1a)-(1c) in $[0,T]$. Then 
\begin{equation*}
\sup_{t \in [0,T]} \left(\lVert w(t)\rVert_{L^{2}}^{2} + \lVert \Lambda^{1+r_{2}}b(t)\rVert_{L^{2}}^{2}\right) + \int_{0}^{T} \lVert \Lambda^{r_{1}} w\rVert_{L^{2}}^{2} + \lVert \Lambda^{1+2r_{2}} b\rVert_{L^{2}}^{2} d\tau \lesssim 1.
\end{equation*}
\end{proposition}

\begin{proof}
We take $L^{2}$-inner products on (21) with $w$ and obtain 
\begin{align}
&\frac{1}{2}\partial_{t} \lVert w\rVert_{L^{2}}^{2} + \lVert \Lambda^{r_{1}} w\rVert_{L^{2}}^{2}\\
\leq& \lVert \Lambda^{r_{2} -1} \nabla \times ((b\cdot\nabla) b) \rVert_{L^{2}} \lVert \Lambda^{1-r_{2}} w\rVert_{L^{2}}\nonumber\\
\lesssim& (\lVert \Lambda^{r_{2}} b\rVert_{L^{\frac{2}{r_{2}}}} \lVert \nabla b\rVert_{L^{\frac{2}{1-r_{2}}}} + \lVert b\rVert_{L^{\infty}}\lVert \Lambda^{r_{2}} \nabla b\rVert_{L^{2}}) \lVert \Lambda^{r_{1}} w\rVert_{L^{2}}\nonumber\\
\lesssim& (\lVert b\rVert_{L^{2}}^{\frac{2r_{2} - 1}{2r_{2}}}\lVert \Lambda^{2r_{2}} b\rVert_{L^{2}}^{\frac{1}{2r_{2}}}\lVert \Lambda^{r_{2}} \nabla b\rVert_{L^{2}} + \lVert b\rVert_{L^{2}}^{\frac{2r_{2} - 1}{2r_{2}}}\lVert \Lambda^{2r_{2}} b\rVert_{L^{2}}^{\frac{1}{2r_{2}}}\lVert \Lambda^{r_{2}} \nabla b\rVert_{L^{2}})\lVert \Lambda^{r_{1}} w\rVert_{L^{2}} \nonumber\\
\leq& \frac{1}{4} \lVert \Lambda^{r_{1}} w\rVert_{L^{2}}^{2} + c(1+ \lVert \Lambda^{2r_{2}} b\rVert_{L^{2}}^{2}) \lVert \Lambda^{r_{2}}\nabla  b\rVert_{L^{2}}^{2}\nonumber
\end{align}
by H$\ddot{o}$lder's inequality, Lemma 2.2, Gagliardo-Nirenberg and Young's inequalities and (2). 

On the other hand, we take $L^{2}$-inner products on (1b) with $\Lambda^{2+2r_{2}}b$ to estimate 
\begin{align*}
&\frac{1}{2}\partial_{t} \lVert \Lambda^{1+r_{2}} b\rVert_{L^{2}}^{2} + \lVert \Lambda^{1+2r_{2}} b\rVert_{L^{2}}^{2}\\
=& -\int [\Lambda^{1+r_{2}} ((u\cdot\nabla) b) - u\cdot\nabla\Lambda^{1+r_{2}} b]\cdot\Lambda^{1+r_{2}} b - \Lambda^{1+r_{2}}((b\cdot\nabla) u) \cdot\Lambda^{1+r_{2}} b\nonumber\\
\lesssim& (\lVert \nabla u\rVert_{L^{\frac{2}{1-r_{1}}}}\lVert \Lambda^{1+r_{2}} b\rVert_{L^{\frac{2}{r_{1}}}} + \lVert \Lambda^{1+r_{2}} u\rVert_{L^{\frac{2}{2r_{2} - 1}}}\lVert \nabla b\rVert_{L^{\frac{1}{1-r_{2}}}}) \lVert \Lambda^{1+r_{2}} b\rVert_{L^{2}}\nonumber\\
&+ (\lVert \Lambda^{1+r_{2}} b\rVert_{L^{\frac{2}{r_{1}}}} \lVert \nabla u\rVert_{L^{\frac{2}{1-r_{1}}}} + \lVert b\rVert_{L^{\infty}}\lVert \Lambda^{2+r_{2}} u\rVert_{L^{2}}) \lVert \Lambda^{1+r_{2}} b\rVert_{L^{2}}\nonumber
\end{align*}
by H$\ddot{o}$lder's inequalities, Lemmas 2.3 and 2.2. We further bound this by 
\begin{align}
&\frac{1}{2}\partial_{t} \lVert \Lambda^{1+r_{2}} b\rVert_{L^{2}}^{2} + \lVert \Lambda^{1+2r_{2}} b\rVert_{L^{2}}^{2}\\
\lesssim& (\lVert \Lambda^{r_{1}}\nabla u\rVert_{L^{2}} \lVert \Lambda^{2+r_{2} - r_{1}}b\rVert_{L^{2}} + \lVert u\rVert_{L^{2}}^{\frac{1}{3+r_{1}}}\lVert \Lambda^{3+r_{1}} u\rVert_{L^{2}}^{\frac{2+r_{1}}{3+r_{1}}}\lVert \Lambda^{2r_{2}} b\rVert_{L^{2}}) \lVert \Lambda^{1+r_{2}} b\rVert_{L^{2}}\nonumber\\
&+ (\lVert \Lambda^{2+r_{2} -r_{1}} b\rVert_{L^{2}} \lVert \Lambda^{r_{1}}\nabla u\rVert_{L^{2}}\nonumber\\
& \hspace{10mm} + \lVert b\rVert_{L^{2}}^{\frac{2r_{2} - 1}{2r_{2}}}\lVert \Lambda^{2r_{2}} b\rVert_{L^{2}}^{\frac{1}{2r_{2}}}\lVert u\rVert_{L^{2}}^{\frac{1+r_{1} - r_{2}}{3+r_{1}}}\lVert \Lambda^{3+r_{1}}u\rVert_{L^{2}}^{\frac{2+r_{2}}{3+r_{1}}})\lVert \Lambda^{1+r_{2}}b\rVert_{L^{2}}\nonumber\\
\lesssim& \left(\lVert \Lambda^{r_{1}}\nabla u\rVert_{L^{2}} \lVert \Lambda^{1+2r_{2}}b\rVert_{L^{2}} + (1+ \lVert \Lambda^{r_{1}} w\rVert_{L^{2}})\lVert \Lambda^{2r_{2}} b\rVert_{L^{2}}\right) \lVert \Lambda^{1+r_{2}} b\rVert_{L^{2}}\nonumber\\
&+ \left(\lVert \Lambda^{1+2r_{2}} b\rVert_{L^{2}} \lVert \Lambda^{r_{1}}\nabla u\rVert_{L^{2}} + (1+ \lVert \Lambda^{2r_{2}} b\rVert_{L^{2}})(1+ \lVert \Lambda^{r_{1}}w\rVert_{L^{2}})\right)\lVert \Lambda^{1+r_{2}}b\rVert_{L^{2}}\nonumber\\
\leq& \frac{1}{4}\left(\lVert \Lambda^{r_{1}} w\rVert_{L^{2}}^{2} + \lVert \Lambda^{1+2r_{2}} b\rVert_{L^{2}}^{2}\right)\nonumber\\
&+ c (1+ \lVert \Lambda^{r_{1}} \nabla u\rVert_{L^{2}}^{2} + \lVert \Lambda^{2r_{2}} b\rVert_{L^{2}}^{2})(1+ \lVert \Lambda^{1+r_{2}} b\rVert_{L^{2}}^{2})\nonumber
\end{align}
by the Sobolev embeddings of $\dot{H}^{r_{1}}(\mathbb{R}^{2}) \hookrightarrow L^{\frac{2}{1-r_{1}}}(\mathbb{R}^{2})$, $\dot{H}^{1-r_{1}}(\mathbb{R}^{2})\hookrightarrow L^{\frac{2}{r_{1}}}(\mathbb{R}^{2})$ and $\dot{H}^{2r_{2} - 1}(\mathbb{R}^{2}) \hookrightarrow L^{\frac{1}{1-r_{2}}}(\mathbb{R}^{2})$, Gagliardo-Nirenberg and Young's inequalities and (2). 

Summing (36) and (37), absorbing dissipative and diffusive terms give 
\begin{align*}
&\partial_{t} (\lVert w\rVert_{L^{2}}^{2}+ \lVert \Lambda^{1+r_{2}} b\rVert_{L^{2}}^{2}) + \lVert \Lambda^{r_{1}} w\rVert_{L^{2}}^{2} + \lVert \Lambda^{1+2r_{2}} b\rVert_{L^{2}}^{2}\\
\lesssim&  (1+ \lVert \Lambda^{r_{1}} \nabla u\rVert_{L^{2}}^{2} + \lVert \Lambda^{2r_{2}} b\rVert_{L^{2}}^{2})(1+ \lVert w\rVert_{L^{2}}^{2} + \lVert \Lambda^{1+r_{2}} b\rVert_{L^{2}}^{2}).
\end{align*}
By (2) and Proposition 4.1, Gronwall's inequality completes the proof of Proposition 4.2. 
\end{proof}

\begin{proposition}
Suppose $\nu, \eta > 0, 0 < r_{1} < \frac{1}{2}, \frac{1}{2} < r_{2} < 1$ so that $r_{1} + r_{2} = 1, $ and $(v, b)$ solves (1a)-(1c) in $[0,T]$. Then 
\begin{align}
\sup_{t \in [0,T]} \lVert \Lambda^{2+r_{2}}b(t)\rVert_{L^{2}}^{2} + \int_{0}^{T} \lVert \Lambda^{2+2r_{2}} b\rVert_{L^{2}}^{2} d\tau \lesssim 1.\nonumber
\end{align}
\end{proposition}

\begin{proof}
We take $L^{2}$-inner products of (1b) with $\Lambda^{4+2r_{2}}$ to estimate 
\begin{align}
&\frac{1}{2} \partial_{t} \lVert \Lambda^{2+r_{2}} b\rVert_{L^{2}}^{2} + \lVert \Lambda^{2+2r_{2}} b\rVert_{L^{2}}^{2}\\
=& -\int [\Lambda^{2+r_{2}}((u\cdot\nabla) b) - u\cdot\nabla\Lambda^{2+r_{2}} b]\cdot\Lambda^{2+r_{2}} b\nonumber\\
& \hspace{10mm} + \int (b\cdot\nabla) u \cdot\Lambda^{4+2r_{2}} b := V_{1} + V_{2}.\nonumber
\end{align}

Firstly, 
\begin{align*}
V_{1} \lesssim& (\lVert \nabla u\rVert_{L^{\infty}} \lVert \Lambda^{2+r_{2}} b\rVert_{L^{2}} + \lVert \Lambda^{2+r_{2}} u\rVert_{L^{\frac{2}{1-r_{1}}}} \lVert \nabla b\rVert_{L^{\frac{2}{1-r_{2}}}} ) \lVert \Lambda^{2+r_{2}} b\rVert_{L^{2}}\\
\lesssim& (\lVert \nabla u\rVert_{L^{2}}^{\frac{1}{2}} \lVert \Lambda^{3} u\rVert_{L^{2}}^{\frac{1}{2}} \lVert \Lambda^{1+r_{2}} b\rVert_{L^{2}}^{\frac{r_{2}}{1+r_{2}}} \lVert \Lambda^{2+2r_{2}} b\rVert_{L^{2}}^{\frac{1}{1+r_{2}}}\\
& \hspace{40mm} + \lVert \Lambda^{3} u\rVert_{L^{2}} \lVert \Lambda^{1+r_{2}} b\rVert_{L^{2}}) \lVert \Lambda^{2+r_{2}} b\rVert_{L^{2}}\nonumber
\end{align*}
by H$\ddot{o}$lder's inequality, Lemma 2.3, Gagliardo-Nirenberg inequalities, the Sobolev embeddings of $\dot{H}^{r_{1}}(\mathbb{R}^{2})\hookrightarrow L^{\frac{2}{1-r_{1}}}(\mathbb{R}^{2})$ and $\dot{H}^{r_{2}}(\mathbb{R}^{2})\hookrightarrow L^{\frac{2}{1-r_{2}}}(\mathbb{R}^{2})$. Using (2) and Proposition 4.2, we can furthermore bound this and use Young's inequalities to obtain 
\begin{equation}
V_{1} \lesssim \left(\lVert \Lambda^{2+2r_{2}} b\rVert_{L^{2}}^{\frac{1}{1+r_{2}}} + 1\right) \lVert \Lambda^{2+r_{2}} b\rVert_{L^{2}} 
\leq \frac{1}{4} \lVert \Lambda^{2+2r_{2}} b\rVert_{L^{2}}^{2} + c(1+ \lVert \Lambda^{2+r_{2}} b\rVert_{L^{2}}^{2}).
\end{equation}

On the other hand, 
\begin{align}
V_{2} \leq& \lVert \Lambda^{2} ((b\cdot\nabla) u)\rVert_{L^{2}} \lVert \Lambda^{2+2r_{2}} b\rVert_{L^{2}}\\
\lesssim& (\lVert \Lambda^{2} b\rVert_{L^{\frac{2}{1-r_{2}}}}\lVert \nabla u\rVert_{L^{\frac{2}{r_{2}}}} + \lVert b\rVert_{L^{\infty}} \lVert \Lambda^{3} u\rVert_{L^{2}})\lVert \Lambda^{2+2r_{2}} b\rVert_{L^{2}}\nonumber\\
\lesssim& (\lVert \Lambda^{2+r_{2}} b\rVert_{L^{2}}\lVert \Lambda^{2-r_{2}} u\rVert_{L^{2}} + \lVert b\rVert_{L^{2}}^{\frac{r_{2}}{1+r_{2}}}\lVert \Lambda^{1+r_{2}} b\rVert_{L^{2}}^{\frac{1}{1+r_{2}}}\lVert \Lambda^{3} u\rVert_{L^{2}}) \lVert \Lambda^{2+2r_{2}} b\rVert_{L^{2}}\nonumber \\
\leq& \frac{1}{4} \lVert \Lambda^{2+2r_{2}} b\rVert_{L^{2}}^{2} + c(\lVert \Lambda^{2+r_{2}} b\rVert_{L^{2}}^{2} \lVert \Lambda^{1+r_{1}} u\rVert_{L^{2}}^{2} + \lVert w\rVert_{L^{2}}^{2})\nonumber\\
\leq& \frac{1}{4} \lVert \Lambda^{2+2r_{2}} b\rVert_{L^{2}}^{2} + + c(1+ \lVert \Lambda^{2+r_{2}} b\rVert_{L^{2}}^{2} ) (1+ \lVert \Lambda^{r_{1}} \nabla u\rVert_{L^{2}}^{2})\nonumber
\end{align}
by H$\ddot{o}$lder's inequality, Lemma 2.2, the Sobolev embeddings of $\dot{H}^{r_{2}}(\mathbb{R}^{2}) \hookrightarrow L^{\frac{2}{1-r_{2}}}(\mathbb{R}^{2})$ and $\dot{H}^{1-r_{2}}(\mathbb{R}^{2})\hookrightarrow L^{\frac{2}{r_{2}}}(\mathbb{R}^{2})$, Gagliardo-Nirenberg and Young's inequalities, (2) and Proposition 4.2. Thus, considering (38), (39) and (40), after absorbing the diffusive term, Gronwall's inequality with (2) completes the proof of Proposition 4.3. 
\end{proof}

Next, we need to attain higher regularity on $w$: 

\begin{proposition}
Suppose $\nu, \eta > 0, 0 < r_{1} < \frac{1}{2}, \frac{1}{2} < r_{2} < 1$ so that $r_{1} + r_{2} = 1, $ and $(v, b)$ solves (1a)-(1c) in $[0,T]$. Then with $ w = \nabla \times v$
\begin{equation*}
\sup_{t \in [0,T]} \lVert \Lambda^{r_{1} + 2r_{2}} w(t)\rVert_{L^{2}}^{2} + \int_{0}^{T} \lVert \Lambda^{2} w\rVert_{L^{2}}^{2} d\tau \lesssim 1.
\end{equation*}
\end{proposition}

\begin{proof}
We take $L^{2}$-inner products on (21) with $\Lambda^{2r_{1} + 4r_{2}}w$ to estimate 
\begin{align*}
&\frac{1}{2}\partial_{t} \lVert \Lambda^{r_{1} + 2r_{2}}w\rVert_{L^{2}}^{2} + \lVert \Lambda^{2} w\rVert_{L^{2}}^{2}\\
=& \int [\Lambda^{r_{1} + 2r_{2}}((u\cdot\nabla) w) - u\cdot\nabla\Lambda^{r_{1} + 2r_{2}} w]\Lambda^{r_{1} + 2r_{2}} w + \int (b\cdot\nabla )j \Lambda^{2r_{1} + 4r_{2}} w\\
\lesssim& (\lVert \nabla u\rVert_{L^{\frac{2}{r_{1}}}} \lVert \Lambda^{r_{1} + 2r_{2}} w\rVert_{L^{2}} + \lVert \Lambda^{r_{1} + 2r_{2}} u\rVert_{L^{2}} \lVert \nabla w\rVert_{L^{\frac{2}{1-r_{2}}}} ) \lVert \Lambda^{r_{1} + 2r_{2}} w\rVert_{L^{\frac{2}{1-r_{1}}}}\\
&+ \lVert \Lambda^{2r_{2}} \nabla \times \text{div} (b\otimes b)\rVert_{L^{2}} \lVert \Lambda^{2} w\rVert_{L^{2}}\\
\lesssim& ( \lVert u\rVert_{L^{2}}^{\frac{r_{1} + 1}{3}}\lVert \Lambda^{3} u\rVert_{L^{2}}^{\frac{2-r_{1}}{3}} \lVert \Lambda^{r_{1} + 2r_{2}} w\rVert_{L^{2}} + \lVert u\rVert_{L^{2}}^{\frac{2-r_{2}}{3}}\lVert \Lambda^{3} u\rVert_{L^{2}}^{\frac{1+r_{2}}{3}}\lVert \Lambda^{1+r_{2}} w\rVert_{L^{2}} ) \lVert \Lambda^{2} w\rVert_{L^{2}}\\
& + \lVert b\rVert_{L^{\infty}} \lVert \Lambda^{2+2r_{2}} b\rVert_{L^{2}} \lVert \Lambda^{2} w\rVert_{L^{2}}\\
\lesssim& \lVert \Lambda^{r_{1} + 2r_{2}} w\rVert_{L^{2}}\lVert \Lambda^{2} w\rVert_{L^{2}} + \lVert b\rVert_{L^{2}}^{\frac{1+r_{2}}{2+r_{2}}}\lVert \Lambda^{2+r_{2}}b\rVert_{L^{2}}^{\frac{1}{2+r_{2}}} \lVert \Lambda^{2+2r_{2}} b\rVert_{L^{2}} \lVert \Lambda^{2} w\rVert_{L^{2}}\\
\leq& \frac{1}{2} \lVert \Lambda^{2} w\rVert_{L^{2}}^{2} + c(1 + \lVert \Lambda^{r_{1} + 2r_{2}}w\rVert_{L^{2}}^{2})(1+ \lVert \Lambda^{2+2r_{2}}b\rVert_{L^{2}}^{2})
\end{align*}
by H$\ddot{o}$lder's inequalities, Lemma 2.3, Gagliardo-Nirenberg inequalities, Propositions 4.2 and 4.3, and Young's inequality. After absorbing the dissipative term, Gronwall's inequality using Proposition 4.3 completes the proof of Proposition 4.4. 
\end{proof}

\begin{proposition}
Suppose $\nu, \eta > 0, 0 < r_{1} < \frac{1}{2}, \frac{1}{2} < r_{2} < 1$ so that $r_{1} + r_{2} = 1, $ and $(v, b)$ solves (1a)-(1c) in $[0,T]$. Then 
\begin{equation*}
\sup_{t \in [0,T]}  \lVert \Lambda^{3} b(t)\rVert_{L^{2}}^{2} + \int_{0}^{T}  \lVert \Lambda^{3+r_{2}} b\rVert_{L^{2}}^{2}d\tau \lesssim 1.
\end{equation*}
\end{proposition}

\begin{proof}
We take $L^{2}$-inner products of (1b) with $\Lambda^{6} b$ to estimate 
\begin{align}
&\frac{1}{2}\partial_{t} \lVert \Lambda^{3} b\rVert_{L^{2}}^{2} + \lVert \Lambda^{3+r_{2}} b\rVert_{L^{2}}^{2}\\
=& -\int [\Lambda^{3} ((u\cdot\nabla) b) - u\cdot\nabla\Lambda^{3} b]\cdot\Lambda^{3} b + \int \Lambda^{3-r_{2}}[(b\cdot\nabla) u]\cdot\Lambda^{3+r_{2}} b\nonumber\\
\lesssim& (\lVert \nabla u\rVert_{L^{\frac{2}{r_{2}}}} \lVert \Lambda^{3} b\rVert_{L^{2}} + \lVert \Lambda^{3} u\rVert_{L^{\frac{2}{r_{2}}}} \lVert \nabla b\rVert_{L^{2}}) \lVert \Lambda^{3} b\rVert_{L^{\frac{2}{1-r_{2}}}}\nonumber\\
&+ \left(\lVert b\rVert_{L^{\infty}} \lVert \Lambda^{4-r_{2}}u\rVert_{L^{2}} + \lVert \Lambda^{3-r_{2}} b\rVert_{L^{\frac{2}{1-r_{2}}}} \lVert \nabla u\rVert_{L^{\frac{2}{r_{2}}}}\right) \lVert \Lambda^{3+r_{2}} b\rVert_{L^{2}}\nonumber\\
\lesssim&(\lVert \nabla u\rVert_{L^{2}}^{r_{2}} \lVert \Delta u\rVert_{L^{2}}^{1-r_{2}}\lVert \Lambda^{3} b\rVert_{L^{2}}\nonumber\\
& \hspace{20mm} + \lVert \Lambda^{3} u\rVert_{L^{2}}^{r_{2}} \lVert \Lambda^{4} u\rVert_{L^{2}}^{1-r_{2}}\lVert b\rVert_{L^{2}}^{\frac{2}{3}} \lVert \Lambda^{3} b\rVert_{L^{2}}^{\frac{1}{3}})\lVert \Lambda^{3+r_{2}} b\rVert_{L^{2}}\nonumber\\
&+ (\lVert b\rVert_{L^{2}}^{\frac{1+r_{2}}{2+r_{2}}}\lVert \Lambda^{2+r_{2}}b\rVert_{L^{2}}^{\frac{1}{2+r_{2}}} \lVert u\rVert_{L^{2}}^{\frac{2r_{2}}{4+r_{2}}}\lVert \Lambda^{4+r_{2}} u\rVert_{L^{2}}^{\frac{4-r_{2}}{4+r_{2}}}\nonumber\\
& \hspace{30mm} + \lVert \Lambda^{3} b\rVert_{L^{2}} \lVert \nabla u\rVert_{L^{2}}^{r_{2}} \lVert \Delta u\rVert_{L^{2}}^{1-r_{2}})\lVert \Lambda^{3+r_{2}} b\rVert_{L^{2}}\nonumber
\end{align}
by H$\ddot{o}$lder's inequalities, Lemmas 2.3 and 2.2, Gagliardo-Nirenberg inequalities and the Sobolev embedding of $\dot{H}^{r_{2}}(\mathbb{R}^{2})\hookrightarrow L^{\frac{2}{1-r_{2}}}(\mathbb{R}^{2})$. We can use (2) and Proposition 4.3 that bounds
\begin{align*}
\sup_{t \in [0,T} \lVert u(t)\rVert_{L^{2}} + \lVert b(t)\rVert_{L^{2}} + \lVert \nabla u(t)\rVert_{L^{2}} + \lVert \Lambda^{2+r_{2}} b(t)\rVert_{L^{2}}^{2} \lesssim 1
\end{align*} 
and Propositions 4.4 that bounds
\begin{align*}
&\sup_{t\in [0,T]} \lVert \Delta u(t)\rVert_{L^{2}}^{2} + \lVert \Lambda^{3} u(t)\rVert_{L^{2}}^{2} + \lVert \Lambda^{4} u(t)\rVert_{L^{2}}^{2} + \lVert \Lambda^{4+r_{2}}u(t)\rVert_{L^{2}}^{2}\\
\lesssim& \sup_{t \in [0,T]} \lVert u(t)\rVert_{L^{2}} + \lVert \Lambda^{r_{2}} \nabla w(t)\rVert_{L^{2}}^{2} \lesssim 1
\end{align*}
to further bound (41) by 
\begin{align*}
&\left(\lVert \Lambda^{3} b\rVert_{L^{2}} +  \lVert \Lambda^{3} b\rVert_{L^{2}}^{\frac{1}{3}}\right)\lVert \Lambda^{3+r_{2}} b\rVert_{L^{2}} + \left(1 + \lVert \Lambda^{3} b\rVert_{L^{2}}  \right) \lVert \Lambda^{3+r_{2}} b\rVert_{L^{2}}\nonumber\\
\leq& \frac{1}{2} \lVert \Lambda^{3+r_{2}} b\rVert_{L^{2}}^{2} + c(1+ \lVert \Lambda^{3} b\rVert_{L^{2}}^{2}).
\end{align*}
due to Young's inequality. Hence after absorbing the diffusive term, Gronwall's inequality completes the proof of Proposition 4.5. 
\end{proof}

We can finally show that the initial regularity is preserved: 
\begin{proposition}
Suppose $\nu, \eta > 0, 0 < r_{1} < \frac{1}{2}, \frac{1}{2} < r_{2} < 1$ so that $r_{1} + r_{2} = 1, $ and $(v, b)$ solves (1a)-(1c) in $[0,T]$. Then  
\begin{align}
\sup_{t \in [0,T]} \lVert \Lambda^{3} v(t)\rVert_{L^{2}}^{2} + \int_{0}^{T} \lVert \Lambda^{3+r_{1}} v\rVert_{L^{2}}^{2} d\tau \lesssim 1.\nonumber
\end{align}
\end{proposition}

\begin{proof}
We take $L^{2}$-inner products of (1a) with $\Lambda^{6}v$ to estimate 
\begin{align}
&\frac{1}{2}\partial_{t} \lVert\Lambda^{3} v\rVert_{L^{2}}^{2} + \lVert \Lambda^{3+r_{1}} v\rVert_{L^{2}}^{2} \\
\lesssim& \lVert \Lambda^{3} ((u\cdot\nabla) v) - u\cdot\nabla \Lambda^{3} v\rVert_{L^{\frac{2}{1+r_{1}}}}\lVert \Lambda^{3} v\rVert_{L^{\frac{2}{1-r_{1}}}}\nonumber\\
&+ \sum_{k=1}^{2}\lVert \Lambda^{3}(v_{k} \nabla u_{k})\rVert_{L^{2}} \lVert \Lambda^{3} v\rVert_{L^{2}} +\lVert \Lambda^{4-r_{1}} (b\otimes b)\rVert_{L^{2}} \lVert \Lambda^{3+r_{1}} v\rVert_{L^{2}}\nonumber\\
\lesssim& (\lVert \nabla u\rVert_{L^{\frac{2}{r_{1}}}} \lVert \Lambda^{3} v\rVert_{L^{2}} + \lVert \Lambda^{3} u\rVert_{L^{2}} \lVert \nabla v\rVert_{L^{\frac{2}{r_{1}}}} ) \lVert \Lambda^{3} v\rVert_{L^{\frac{2}{1-r_{1}}}}\nonumber\\
&+ (\lVert \Lambda^{3} v\rVert_{L^{\frac{2}{1-r_{1}}}}\lVert \nabla u \rVert_{L^{\frac{2}{r_{1}}}} + \lVert v\rVert_{L^{\frac{2}{r_{1}}}} \lVert  \Lambda^{3} \nabla u\rVert_{L^{\frac{2}{1-r_{1}}}}) \lVert \Lambda^{3} v\rVert_{L^{2}}\nonumber\\
&+ \lVert b\rVert_{L^{\infty}} \lVert \Lambda^{3+r_{2}} b\rVert_{L^{2}} \lVert \Lambda^{3+r_{1}} v\rVert_{L^{2}}\nonumber
\end{align}
by H$\ddot{o}$lder's inequalities, Lemmas 2.3 and 2.2. We further bound (42) by 
\begin{align*}
&(\lVert \nabla u\rVert_{L^{2}}^{\frac{1+r_{1}}{2}}\lVert \Lambda^{3} u\rVert_{L^{2}}^{\frac{1-r_{1}}{2}} \lVert \Lambda^{3} v\rVert_{L^{2}}\\
& \hspace{20mm} + \lVert w\rVert_{L^{2}} \lVert \nabla v\rVert_{L^{2}}^{\frac{1}{1+r_{2}}}\lVert \Lambda^{1+r_{2}} \nabla v\rVert_{L^{2}}^{\frac{r_{2}}{1+r_{2}}}) \lVert \Lambda^{3+r_{1}} v\rVert_{L^{2}}\\
+&  (\lVert \Lambda^{3+r_{1}} v\rVert_{L^{2}}\lVert \nabla u\rVert_{L^{2}}^{\frac{1+r_{1}}{2}}\lVert \Lambda^{3} u\rVert_{L^{2}}^{\frac{1-r_{1}}{2}}\\
& \hspace{20mm} + (\lVert u\rVert_{L^{2}} + \lVert \Delta u\rVert_{L^{2}})^{r_{1}}\lVert \nabla v\rVert_{L^{2}}^{1-r_{1}} \lVert \Lambda^{3+r_{1}} v\rVert_{L^{2}} )\lVert \Lambda^{3} v\rVert_{L^{2}}\\
+&  (\lVert b\rVert_{L^{2}} + \lVert \Lambda^{3} b\rVert_{L^{2}}) \lVert \Lambda^{3+r_{2}} b\rVert_{L^{2}} \lVert \Lambda^{3+r_{1}} v\rVert_{L^{2}}\\
\lesssim& \left(\lVert \Lambda^{3} v\rVert_{L^{2}} + \lVert w\rVert_{L^{2}}^{\frac{1}{1+r_{2}}}\lVert \Lambda^{1+r_{2}} w\rVert_{L^{2}}^{\frac{r_{2}}{1+r_{2}}}\right) \lVert \Lambda^{3+r_{1}} v\rVert_{L^{2}}\\
+&  \left(\lVert \Lambda^{3+r_{1}} v\rVert_{L^{2}}+ (1 + \lVert w\rVert_{L^{2}})^{r_{1}}\lVert w\rVert_{L^{2}}^{1-r_{1}} \lVert \Lambda^{3+r_{1}} v\rVert_{L^{2}} \right)\lVert \Lambda^{3} v\rVert_{L^{2}}\\
+& \lVert \Lambda^{3+r_{2}} b\rVert_{L^{2}} \lVert \Lambda^{3+r_{1}} v\rVert_{L^{2}}\\
\lesssim& \left(\lVert \Lambda^{3} v\rVert_{L^{2}} + 1\right) \lVert \Lambda^{3+r_{1}} v\rVert_{L^{2}}\\
& \hspace{20mm} + \lVert \Lambda^{3+r_{1}} v\rVert_{L^{2}}\lVert \Lambda^{3} v\rVert_{L^{2}} + \lVert \Lambda^{3+r_{2}} b\rVert_{L^{2}} \lVert \Lambda^{3+r_{1}} v\rVert_{L^{2}}\\
\leq& \frac{1}{2} \lVert \Lambda^{3+r_{1}} v\rVert_{L^{2}}^{2} + c(1+ \lVert \Lambda^{3} v\rVert_{L^{2}}^{2}) (1+ \lVert \Lambda^{3+r_{2}} b\rVert_{L^{2}}^{2})
\end{align*}
by Gagliardo-Nirenberg inequalities, the Sobolev embedding of $H^{3}(\mathbb{R}^{2})\hookrightarrow L^{\infty}(\mathbb{R}^{2})$ and $\dot{H}^{r_{1}} (\mathbb{R}^{2}) \hookrightarrow L^{\frac{2}{1-r_{1}}}(\mathbb{R}^{2})$, Lemma 2.1, (2), Propositions 4.2, 4.3 and 4.4 and Young's inequality. After absorbing the dissipative term, Gronwall's inequality using Proposition 4.5 completes the proof of Proposition 4.6. 
\end{proof}

\section{Proof of Theorem 1.1}

With \textit{a priori } estimates achieved in previous sections, specifically Propositions 3.3, 3.5, 4.5 and 4.6, it is a standard procedure to complete the proof of Theorem 1.1. We recall the mollification of $\mathcal{J}_{\epsilon}f$ of $f \in L^{p}(\mathbb{R}^{2}), 1 \leq p \leq \infty$ by 
\begin{equation*}
(\mathcal{J}_{\epsilon}f)(x) = \epsilon^{-2} \int_{\mathbb{R}^{2}} \rho\left(\frac{x-y}{\epsilon}\right)f(y) dy, \hspace{3mm} \epsilon > 0
\end{equation*}
where $\rho (\lvert x\rvert) \in C_{0}^{\infty}, \rho \geq 0, \int_{\mathbb{R}^{2}} \rho dx = 1$. We regularize (1a)-(1c) as follows:
\begin{align*}
\begin{cases}
\partial_{t} v^{\epsilon} + \mathcal{J}_{\epsilon}((\mathcal{J}_{\epsilon}u^{\epsilon}) \cdot \nabla (\mathcal{J}_{\epsilon} v^{\epsilon})) + \sum_{k=1}^{2} \mathcal{J}_{\epsilon}((\mathcal{J}_{\epsilon} v_{k}^{\epsilon} \nabla (\mathcal{J}_{\epsilon} u_{k}^{\epsilon}))\\
\hspace{5mm} + \nabla (\pi^{\epsilon} + \frac{1}{2} \lvert b^{\epsilon}\rvert^{2}) + \mathcal{J}_{\epsilon} (\Lambda^{2r_{1}} \mathcal{J}_{\epsilon} v^{\epsilon}) = \mathcal{J}_{\epsilon} ((\mathcal{J}_{\epsilon} b^{\epsilon})\cdot\nabla) \mathcal{J}_{\epsilon} b^{\epsilon})),\\
\partial_{t}b^{\epsilon} + \mathcal{J}_{\epsilon} ((\mathcal{J}_{\epsilon} u^{\epsilon})\cdot\nabla( \mathcal{J}_{\epsilon} b^{\epsilon})) - \mathcal{J}_{\epsilon} ((\mathcal{J}_{\epsilon}b^{\epsilon})\cdot\nabla( \mathcal{J}_{\epsilon}u^{\epsilon})) + \mathcal{J}_{\epsilon} (\Lambda^{2r_{2}}\mathcal{J}_{\epsilon} b^{\epsilon}) = 0,\\
\nabla\cdot u^{\epsilon} = \nabla \cdot b^{\epsilon} = 0, \hspace{3mm} v^{\epsilon} = (1-\Delta) u^{\epsilon},
\end{cases}
\end{align*}
so that using properties of mollifiers, one can show via Picard Theorem, the global existence of the regularized solution pair which will lead through the process of obtaining a uniform bound locally in time and then using Alaoglu's theorem, the existence of the local solution pair to (1a)-(1c). We omit further details referring to [21] and [23].


\begin{thebibliography}{100}  
\addtolength{\leftmargin}{0.2in} 
\setlength{\itemindent}{-0.2in} 

\bibitem[1]{1} C. Cao, and J. Wu, \textit{Global regularity for the 2D MHD equations with mixed partial dissipation and magnetic diffusion}, Adv. Math., \textbf{226} (2011), 1803-1822. 

\bibitem[2]{2} C. Cao, J. Wu, and B. Yuan, \textit{The 2D incompressible magnetohydrodynamics equations with only magnetic diffusion}, arXiv:1306.3629 [math.AP], 21, Jan., 2014. 

\bibitem[3]{3} D. Chae, \textit{Global regularity for the 2-D Boussinesq equations with partial viscous terms}, Adv. Math., \textbf{203}, 2 (2006), 497-513.

\bibitem[4]{4} J.-Y. Chemin, \textit{Perfect incompressible fluids}, Oxford lecture series in mathematics and its applications, \textbf{14}, Oxford University Press Inc., New York (1998).

\bibitem[5]{5} S. Chen, C. Foias, D. D. Holm, E. Olson, E. S. Titi, and S. Wynne, \textit{The Camassa-Holm equations as a closure model for turbulent channel flow}, Phys. Rev. Lett., \textbf{81}, 24 (1998), 5338-5341. 

\bibitem[6]{6} A. Cheskidov, D. D. Holm, E. Olson, and E. S. Titi, \textit{On a Leray$-alpha$ model of turbulence}, Proc. R. Soc. Lond. Ser. A Math. Phys. Eng. Sci., \textbf{461} (2005), 629-649. 

\bibitem[7]{7} P. Constantin, G. Iyer, and J. Wu,  \textit{Global regularity for a modified critical dissipative quasi-geostrophic equation}, Indiana Univ. Math. J., \textbf{57} (2001), 97-107, Special Issue.

\bibitem[8]{8} A. C$\acute{o}$rdoba, and D. C$\acute{o}$rdoba,  \textit{A maximum principle applied to quasi-geostrophic equations}, Comm. Math. Phys., \textbf{249}, 3 (2004), 511-528.

\bibitem[9]{9} J. Fan, and T. Ozawa, \textit{Global Cauchy problem for the 2-D magnetohydrodynamic-$\alpha$ models with partial viscous terms}, J. Math. Fluid Mech., \textbf{12} (2010), 306-319.

\bibitem[10]{10} J. D. Gibbon, and D. D. Holm, \textit{Estimates for the LANS$-\alpha$, Leray$-alpha$ and Bardina models in terms of a Navier-Stokes Reynolds number}, Indiana Univ. Math. J., \textbf{57}, 6 (2008), 2761-2773, Special Issue. 

\bibitem[11]{11} T. Hmidi, and S. Keraani, \textit{On the global well-posedness of the Boussinesq system with zero viscosity}, Indiana Univ. Math. J., \textbf{58}, 4 (2009), 1591-1618.

\bibitem[12]{12} T. Hou, and C. Li, \textit{Global well-posedness of the viscous Boussinesq equations}, Discrete Contin. Dyn. Syst., \textbf{12}, 1 (2005), 1-12.

\bibitem[13]{13} A. A. Ilyin, E. M. Lunasin, and E. S. Titi, \textit{A modified-Leray-$\alpha$ subgrid scale model of turbulence}, Nonlinearity, \textbf{29} (2006), 879-897.  

\bibitem[14]{14} Q. Jiu, and J. Zhao, \textit{Global regularity of 2D generalized MHD equations with magnetic diffusion}, arXiv:1309.5819 [math.AP], 21, Jan., 2014. 

\bibitem[15]{15} Q. Jiu, and J. Zhao, \textit{A remark on global regularity of 2D generalized magnetohydrodynamic equations}, arXiv:1306.2823 [math.AP], 21, Jan., 2014. 

\bibitem[16]{16} N. Ju, \textit{The maximum principle and the global attractor for the dissipative 2D quasi-geostrophic equations}, Comm. Math. Phys., \textbf{255}, 1 (2005), 161-181.

\bibitem[17]{17} T. Kato, \textit{Liapunov functions and monotonicity in the Navier-Stokes equation}, Functional-analytic methods for partial differential equations, lecture notes in mathematics, \textbf{1450} (1990), 53-63. 

\bibitem[18]{18} T. Kato, and G. Ponce, \textit{Commutator estimates and the Euler and Navier-Stokes equations}, Comm. Pure Appl. Math., \textbf{41}, 7 (1988), 891-907.

\bibitem[19]{19} A. Kiselev, \textit{Regularity and blow up for active scalars}, Math. Model. Nat. Phenom., \textbf{5}, 4 (2010), 225-255. 

\bibitem[20]{20} J. S. Linshiz and E. S. Titi, \textit{Analytical study of certain magnetohydrodynamic-$\alpha$ models}, J. Math, Phys., \textbf{48}, 065504 (2007).

\bibitem[21]{21} A. J. Majda, and A. L. Bertozzi, \textit{Vorticity and incompressible flow}, Cambridge University Press, Cambridge, 2001.

\bibitem[22]{22} E. Olson, and E. S. Titi, \textit{Viscosity versus vorticity stretching: global well-posedness for a family of Navier-Stokes-alpha-like models}, Nonlinear Anal., \textbf{66} (2007), 2427-2458. 

\bibitem[23]{23} M. Sermange, and R. Temam, \textit{Some mathematical questions related to the MHD equations}, Comm. Pure Appl. Math., \textbf{36} (1983), 635-664.  
  
\bibitem[24]{24} C. V. Tran, X. Yu, and Z. Zhai, \textit{On global regularity of 2D generalized magnetohydrodynamics equations}, J. Differential Equations, \textbf{254}, 10 (2013), 4194-4216.

\bibitem[25]{25} C. V. Tran, X. Yu, and Z. Zhai, \textit{Note on solution regularity of the generalized magnetohydrodynamic equations with partial dissipation}, Nonliner Anal., \textbf{85} (2013), 43-51. 

\bibitem[26]{26} J. Wu, \textit{The generalized MHD equations}, J. Differential Equations, \textbf{195} (2003), 284-312. 

\bibitem[27]{27} J. Wu, \textit{Global regularity for a class of generalized magnetohydrodynamic equations}, J. Math. Fluid Mech., \textbf{13}, 2 (2011), 295-305. 

\bibitem[28]{28} K. Yamazaki, \textit{Remarks on the global regularity of two-dimensional magnetohydrodynamics system with zero dissipation}, Nonlinear Anal., \textbf{94} (2014), 194-205.   

\bibitem[29]{29} K. Yamazaki, \textit{Regularity criteria of supercritical beta-generalized quasi-geostrophic equation in terms of partial derivatives}, Electron. J. Differential Equations, \textbf{2013}, 217 (2014), 1-12. 

\bibitem[30]{30} K. Yamazaki, \textit{Global regularity of logarithmically supercritical MHD system with zero diffusivity}, Appl. Math. Lett., \textbf{29} (2014), 46-51. 

\bibitem[31]{31} K. Yamazaki, \textit{On the global regularity of two-dimensional generalized magnetohydrodynamics system}, arXiv:1306.2842 [math.AP], 21, Jan., 2014.

\bibitem[32]{32} B. Yuan, and L. Bai, \textit{Remarks on global regularity of 2D generalized MHD equations}, arXiv:1306.2190 [math.AP], 21, Jan., 2014.

\bibitem[33]{33} J. Zhao, and M. Zhu, \textit{Global regularity for the incompressible MHD$-\alpha$ system with fractional diffusion}, Appl. Math. Lett., \textbf{29} (2014), 26-29. 

\bibitem[34]{34} Y. Zhou, and J. Fan, \textit{Regularity criteria for the viscous Camassa-Holm equations}, Int. Math. Res. Not. IMRN, 2009 (2009), 2508-2518. 

\bibitem[35]{35} Y. Zhou, and J. Fan, \textit{On the Cauchy problem for a Leray$-\alpha$-MHD model}, Nonlinear Anal. Real World Appl., \textbf{12} (2011), 648-657.

\bibitem[36]{36} Y. Zhou, and J. Fan, \textit{Regularity criteria for a magnetohydrodynamic$-\alpha$ model}, Commun. Pure Appl. Anal., \textbf{10}, 1 (2011), 309-326. 

\end{thebibliography}
\end{document}